\documentclass[10pt,twocolumn,doublespace]{IEEEtran}

\usepackage{amsmath,amssymb,epsfig,psfrag,cite,subfigure}

\include{macros}
\usepackage{bbm}

\usepackage{graphicx}
\usepackage{algorithmic}
\usepackage{algorithm}
\usepackage{epsfig,psfrag}
\usepackage{subfigure}
\usepackage{color}
\usepackage{url}
\newtheorem{remark}{Remark}

\usepackage{pstricks, pst-node, pst-plot, pst-circ}
\usepackage{moredefs}
\usepackage{tikz}
\tikzstyle{superpeers}=[draw,circle,burntorange, left color=\oran,
                       text=violet,minimum width=10pt]

\newtheorem{assumption}{Assumption}

\graphicspath{{figures/}}

\providelength{\AxesLineWidth}       
\providelength{\GridLineWidth}       
\providelength{\GridLineDotSep}      
\providelength{\MinorGridLineWidth}  
\providelength{\MinorGridLineDotSep} 
\providelength{\plotwidth}           
\providelength{\LineWidth}           
\providelength{\MarkerSize}          
\newrgbcolor{GridColor}{0.8 0.8 0.8}%

\begin{document}

\title{Diffusion Estimation Over Cooperative Multi-Agent Networks With Missing Data}
\author{{Mohammad  Reza Gholami,  \emph{Member, IEEE}, Magnus Jansson, \emph{Member, IEEE},\\
 Erik G. Str\"om, \emph{Senior Member, IEEE}, and Ali H. Sayed, \emph{Fellow, IEEE}}

\thanks{M. R. Gholami is with Campanja AB, Stockholm SE-111 57, Sweden
(e-mail: reza@campanja.com)}
\thanks{M. Jansson is with the ACCESS Linnaeus Center, Electrical Engineering, KTH--Royal Institute of Technology, SE-100 44 Stockholm, Sweden
(e-mail: janssonm@kth.se)}
\thanks{E. G. Str\"om is with the Division of Communication Systems, Department of Signals and Systems, Chalmers University of Technology, SE-412 96 Gothenburg, Sweden (e-mail:erik.strom@chalmers.se).}
\thanks{A. H. Sayed is with the Department of Electrical Engineering, University
of California, Los Angeles, CA 90095 (e-mail: sayed@ee.ucla.edu).}
 \thanks{This work was initiated while M. R. Gholami was a visiting student at the UCLA Adaptive Systems Laboratory. The
work was supported in part by the Swedish Research Council (contract no.~2007-6363) and the Access center
and in part by Ericsson's Research Foundation
(Reference: FOSTIFT-12:021 and FOSTIFT-13:031). The work of A. H. Sayed was supported in part by NSF grants CCF-1524250 and 
ECCS-1407712. A short and preliminary version of this work was presented in the conference publication~\cite{Gholami_diff_globalsip_2013}.
 }
}

\maketitle

\begin{abstract}
In many fields, and especially in the medical and social sciences and in recommender systems, data are gathered through clinical studies or targeted surveys. Participants are generally reluctant to respond to all questions in a survey or they may lack information to respond adequately to some questions. The data collected from these studies tend to lead to linear regression models where the regression vectors are only known partially: some of their entries are either missing completely or replaced randomly by noisy values.
In this work, assuming missing positions are replaced by noisy values, we examine how a connected network of agents, with each one of them subjected to a stream of data with incomplete regression information, can cooperate with each other through local interactions to estimate the underlying model parameters in the presence of missing data.
We explain how to adjust the distributed diffusion through (de)regularization in order to eliminate the bias introduced by the incomplete model.
We also propose a technique to recursively estimate the (de)regularization parameter
and examine the performance of the
resulting strategy. We illustrate the results by considering two applications: one dealing with a mental health survey and
the other dealing with a household consumption survey.


\end{abstract}
\section{Introduction}

In data gathering procedures through clinical studies or targeted surveys, it is common that some components of the data are missing or left unobserved.
For example, in a clinical study, a participant may be reluctant to answer some questions or may drop out
of the survey and never return it~\cite{enders2010applied}. Likewise, in a recommender system using content based filtering~\cite{koren2009matrix}, a participant may
prefer to leave some questions unanswered.
The phenomenon of missing data is prevalent in many
fields including the social sciences, medical sciences, econometrics, survival data analysis, and machine learning \cite{ibrahim2005missing,bernhardt2013statistical,rigobon2007estimation,cox1984analysis,enders2006primer,enders2010applied,marlin2008missing,baraldi2010introduction,robins1994estimation}.
In general, three patterns are considered for missing data: missing at random (MAR), missing completely at random (MCAR), and missing not at random (MNAR)~\cite{davey2009statistical,little2002statistical,rigobon2009bias}.
In MAR patterns, the probability of missing a variable in data gathering is related not only to the value of the variable but also to some other measured variables~\cite{enders2010applied}.
In comparison, in MCAR scenarios, the probability of missing a variable is solely related to the value of that variable~\cite{davey2009statistical}.
For example, if an observation is below a threshold, it may not be observed.
In MNAR, the probability of missing depends on unobserved values. For example, in a cancer trial
some patients may be feeling weak and unable to participate
in the study~\cite{enders2010applied}. In general it may not be easy to verify which pattern of missing is happening for a study.
For a survey on different missing patterns, the reader may refer to, e.g., \cite[Ch. 1]{enders2010applied} and the references therein.

There have been many studies on inference problems for
missing data with several useful techniques proposed to address the challenges associated with censored data.
Many of the approaches are based on heuristic methods, as already noted in~\cite{loh2012corrupted}. There are
broadly two classes of techniques to deal with estimation in the presence of missing data: imputation and deletion (listwise or pairwise)~\cite{arminger1995handbook,wang2006efficiencies,wang2012multiple,royston2004multiple,enders2010applied,rubin2009multiple,baraldi2010introduction}.
If the positions of the missing data are known in advance,
then they can either be replaced by some deterministic or random values (sometimes called single or multiple imputation), or the corresponding data can be removed from the dataset altogether. Removing data generally leads to performance degradation,
although this practice is accepted in some statistical software packages,
e.g., SAS, SPSS, GAUSS, and CDAS, mainly
because of ease of implementation. In contrast, data imputation results in biased estimates~\cite{arminger1995handbook,rigobon2007estimation,rigobon2009bias,royston2004multiple}.
For a discussion on different types of deletion and imputation in missing data analysis, see \cite[Ch. 2]{enders2010applied}.

When data imputation is assumed, one may resort to a maximum likelihood procedure to estimate the missing data
if the distribution of the data happens to be known beforehand~\cite{horton1999maximum,trawinski1964maximum}.
For example, the expectation-maximization (EM) technique  provides one useful solution method~\cite{marlin2008missing}.
However, the EM implementation tends to be computationally intensive and can suffer from convergence issues~\cite{wiesel2008linear}.
If the positions of the missing data are unknown, a mixture model may be used to describe the data model \cite{bishop:2006:PRML} and the EM algorithm can be subsequently applied to estimate the missing data \cite{EM_mixture,ghahramani1994supervised}. However, the number of components can become large in high dimensions.
For other possible recovery techniques including ordinary least-squares, inverse probability weighting, and fully Bayseian methods, the reader may refer to, e.g.,
\cite{li2011weighting,ibrahim2005missing,tsimikas2012inference,nelson1972linear,nelson1973linear}.

In this work, we examine how a connected network of agents, with each one of them subjected to a
stream of data with incomplete regression information, can cooperate with each other to estimate the
underlying model parameters in the presence of missing data.
In particular, we assume that some entries of the regressor can be missed according to the MCAR pattern,
where the missing data are replaced by random entries.
Note that for nonstreaming data, there are various techniques in the literature \cite{consensus_2013,dogandzic2008decentralized,fang2010optimal} to handle an inference problem.
These techniques are useful if agents
are able to collect the data and save them in buffers for batch or centralized processing.

Indeed, one of the main approaches used in the literature to address estimation in the presence of missing data relies on centralized solutions, where the data is collected centrally and processed by a fusion center. This approach has the
disadvantage of requiring a large amount of energy and communication and may limit scalability and robustness of the network.
There are many situations, however, where the data is already available in a distributed manner (e.g., separate clinics collect surveys from their patients independently of each other). Moreover, in many of these cases, privacy and secrecy considerations may prevent the
clinics or survey-collectors to share or transmit the data to a central location. In these situations,
it is  preferable to process the data in a distributed manner.
In this approach, only estimates of the unknown parameter model, and not the raw collected data, are exchanged between neighboring nodes.
This paper focuses on this approach and develops a distributed strategy for handling missing data from results collected at dispersed locations. One of the key challenges is to infuse the distributed procedure with a mechanism to combat the bias that results from the presence of censored data.

For this purpose, we shall rely on the idea of \textit{diffusion} adaptation~\cite{sayed_2014,Sayed_2014_adptive,sayed_2013_magazine,sayed_2012_diff_net,Jianshu_2012}
to design a fully distributed algorithm
that is able to adapt and learn from streaming data.
Useful consensus strategies can also be employed \cite{consensus_2013,Kar_Mura_2011,Dimakis_2010,Lorenzo_2011,Braca_2008,fang2010distributed}. We focus on diffusion strategies in this work due to their proven enhanced stability ranges and improved mean-square-performance properties over consensus  networks when constant step-sizes are employed for continuous adaptation and learning \cite{tu2012diffusion,sayed_2013_magazine,Sayed_2014_adptive,sayed_2012_diff_net}; consensus networks are sensitive to the network topology and their state error vectors can grow unbounded when constant step-sizes are used. We will explain how to adjust the diffusion strategy through (de)regularization in order to eliminate
the bias introduced by imputation.
We will examine the stability and performance of the proposed algorithm
and provide computer simulations on two applications: one dealing with a mental health survey and the other dealing with a household consumption survey.

\textbf{Notation.} We use lowercase letters to denote vectors, uppercase letters for matrices, plain letters for deterministic variables, boldface letters for random variables. We use $\odot$ and $\otimes$ for the Hadamard (elementwise) and Kronecker products, respectively.
In addition, $\mathrm{diag}\{x_1,\ldots,x_N\}$ denotes a diagonal matrix with diagonal elements $x_1,\ldots,x_N$.
We use $\mathrm{col}\{a,b\}$ to represent a column vector with entries $a$ and $b$, while
$I_M$ and $0_M$ denote the $M \times M$ identity and $M \times M$ zero matrices.
We use $\lambda_i(A)$ to denote the $i-$th eigenvalue of matrix $A$.
The $i-$th and $j-$th entry of a matrix $A$ is represented by $A(i,j)$.

\section{Problem statement}
\label{sec:model}
Consider a connected network with $N$ agents. Each agent senses a stream of wide-sense stationary data  $\{\boldsymbol{d}_k(i),\boldsymbol{u}_{k,i}\}$ that satisfy the linear regression model:
\begin{align}
\label{eq:data_model}
\boldsymbol{d}_k(i)=\boldsymbol{u}_{k,i} w^o+\boldsymbol{v}_k(i), \quad k=1,\ldots, N
\end{align}
where $k$ is the node index and $i$ is the time index. The row vector $\boldsymbol{u}_{k,i}$
denotes a zero-mean random process
with covariance matrix $R_{u,k}=\mathbb{E}\,\boldsymbol{u}_{k,i}^*\boldsymbol{u}_{k,i}>0$, while
$\boldsymbol{v}_k(i)$ is a zero-mean white noise process with variance $\sigma_{v,k}^2$. The column vector
$w^o\in \mathbb{R}^M$ is
the unknown parameter that the network is interested in estimating.
\begin{remark}
\label{re:rem1}
Sometimes, as illustrated by the examples discussed later in Sec. V, we may encounter data models of the form
\begin{align}
\label{eq:data_model_prac}
\bar{\boldsymbol{d}}_k(i)=\alpha+\bar{\boldsymbol{u}}_{k,i} w^o+\boldsymbol{v}_k(i), \quad k=1,\ldots, N
\end{align}
where $\alpha$ is some deterministic intercept value, and $\bar{\boldsymbol{u}}_{k,i}$ is a row vector with nonzero mean.
This model can still be reduced to \eqref{eq:data_model} by subtracting the mean of the measurement from both sides of \eqref{eq:data_model_prac}.
\end{remark}

\begin{assumption}
\label{Assum:1}
Continuing with model \eqref{eq:data_model}, we assume  that the regression and noise processes are each spatially independent and temporally white.
In addition,
we assume that $\boldsymbol{u}_{\ell,i}$ and $\boldsymbol{v}_k(j)$ are independent
of each other for all $\ell, i, k$, and $j$.$\qquad \qquad \qquad \qquad \qquad\qquad\qquad \qquad\qquad \qquad \qquad\qquad  \blacksquare$%
\end{assumption}

In this study, we examine the situation in which some entries in the regression vectors may be missing
completely at random due to a variety of reasons, including incomplete information or censoring. We
denote the incomplete regressor by $\bar{\boldsymbol{u}}_{k,i}$ and express it in the form:
\begin{align}
\label{eq:relation_data}
\bar{\boldsymbol{u}}_{k,i}=\boldsymbol{u}_{k,i}(I_M-\boldsymbol{F}_{k,i})+\boldsymbol{\xi}_{k,i} \boldsymbol{F}_{k,i}
\end{align}
where $\boldsymbol{F}_{k,i}=\mathrm{diag}\big\{\boldsymbol{f}^1_{k,i},\ldots,\boldsymbol{f}^M_{k,i}\big\}$ consists of random indicator variables,
 $\boldsymbol{f}^j_{k,i}\in\{0,1\}$.
Each variable $\boldsymbol{f}_{k,i}^j$ is equal to one with some probability $p<1$ and equal to zero with probability $1-p$. The value of $p$ represents the likelihood that the $j-$th entry of the regression vector $\boldsymbol{u}_{k,i}$ is missing at time $i$. In that case, the missing entry is assumed to be replaced by an entry from the zero-mean perturbation (row) variable~$\boldsymbol{\xi}_{k,i}$.

The model considered in \eqref{eq:relation_data} encompasses
different types of censoring such as \emph{left censoring}, \emph{right censoring}, \emph{interval censoring}, and
\emph{random censoring}~\cite{efron_censored_1981,censored_1976,censored_1979}.\footnote{In econometrics,
the use of the ``\emph{coded}'' terminology (such as top-coded and bottom-coded) is more common than censoring \cite{rigobon2007estimation,chernozhukov2009set}.}
In the missing data literature, the position of missing data is often assumed to be
known in advance \cite{arminger1995handbook}.
In general, both scenarios with known and unknown missing positions may happen in practice. For example,
a participant in a survey may leave a question unanswered (known missing position), or may report a wrong value
(unknown missing position since the true value is intentionally replaced by a wrong value); the latter may be considered an outlier although usually an outlier is
defined as a measurement that is distant from other measurements.
For a known missing position, the imputation technique can then be used to fill the position of data left unobserved yielding
a modified data model.
If the analyst for the modified model is not aware of the imputed position, we can still consider
model \eqref{eq:relation_data} with the assumption that the imputer and the analyst are different and do not exchange information about the missing position.

\begin{remark}
Model \eqref{eq:relation_data} is sometimes referred to as a mixture model in the literature  \cite{hedeker1997application,bailey1994fitting,little1993pattern,ghahramani1994supervised}. 
In such models, it is mainly assumed that the  components follow Gaussian distributions with different means or variances (or both). However, in this study, we have no particular assumptions about the distribution of the components or their first and second-order moments. For mixture model approaches, if the distribution of the data deviates from the normal distribution, a large number of Gaussian components is then mixed to model the missing data \cite{rasmussen1999infinite,diebolt1994estimation}. In these cases, the design of the algorithms can become more challenging.
\end{remark}

\begin{assumption}
\label{Assum:2}
We assume that the random variables $\boldsymbol{u}_{k,i}, \boldsymbol{f}^j_{k,i}$, and $\boldsymbol{\xi}_{k,i}$ are independent of each other. We also assume that the random process $\boldsymbol{\xi}_{k,i}$ is temporally white and spatially independent with covariance matrix ${\mathbb{E}\, \boldsymbol{\xi}_{k,i}^*\boldsymbol{\xi}_{k,i}=\sigma_{\xi,k}^2I_M}$.~$\qquad \qquad \qquad \qquad \qquad \qquad \qquad \quad~~\blacksquare$
\end{assumption}

From model \eqref{eq:data_model}, the minimum mean-square-error (MMSE) estimate of the unknown vector $w^o$ based on the data collected at node $k$ is given
by~\cite{kailath2000linear}:
\begin{align}
\label{eq:opt_est}
w_k^o=R^{-1}_{u,k}\,r_{du,k}
\end{align}
where
\begin{align}
{r_{du,k}\triangleq \mathbb{E}\,\boldsymbol{d}_k(i)\boldsymbol{u}^*_{k,i}}.
\end{align}
It is easy to verify from \eqref{eq:data_model} that ${w_{k}^o=w^o}$ so that the MMSE solution allows node $k$ to recover the unknown $w^o$ if the actual moments $\{R_{u,k}, r_{du,k}\}$ happen to be known.  The resulting mean-square-error is~\cite{kailath2000linear}:
\begin{align}
\label{eq:min_mse_perfect}
J_{k,{\mathrm{min}}}\triangleq J_k(w^o_k)&=\mathbb{E}|\boldsymbol{d}_{k}(i)-\boldsymbol{u}_{k,i}w^o_k|^2\nonumber\\
&=\sigma^2_{d,k}-r^*_{du,k}R_{u,k}^{-1}r_{du,k}\nonumber\\
&=\sigma^2_{v,k}.
\end{align}

Let us now investigate the estimate that would result
if we use instead the censored regressor,
$\bar{\boldsymbol{u}}_{k,i}$.
The estimate in this case is given by
\begin{align}
\label{eq:op_missing}
 \bar{w}_k^o=R^{-1}_{\bar{u},{k}}\,r_{d\bar{u},{k}}
\end{align}
with the covariance matrix $R_{\bar{u},k}$ computed as
\begin{align}
\label{eq:cov_censored}
R_{\bar{u},k}&\triangleq\mathbb{E}\bar{\boldsymbol{u}}^*_{k,i}\bar{\boldsymbol{u}}_{k,i}\nonumber\\
&\mathop{=}\limits^{(a)}
\mathbb{E}\{(I-\boldsymbol{F}_{k,i})R_{u,k}(I-\boldsymbol{F}_{k,i})\}+\sigma^2_{\xi,k}\mathbb{E}\{\boldsymbol{F}^2_{k,i}\}\nonumber\\
&=R_{u,k}-P_1 \odot R_{u,k}+p\sigma^2_{\xi,k}I_M\nonumber\\
&\triangleq R_{u,k}+R_{r,k}
\end{align}
where in $(a)$ we used the independence of $\boldsymbol{u}_{k,i}$ and $\boldsymbol{\xi}_{k,i}$ (Assumption \ref{Assum:2})
and where the matrix $R_{r,k}$ is defined as
\begin{align}
\label{eq:R_2}
R_{r,k}\triangleq-P_1 \odot R_{u,k}+p\sigma_{\xi,k}^2I_M
\end{align}
in terms of
\begin{eqnarray}
P_{1}&\triangleq& (2p-p^2)\mathbbm{1}_M\mathbbm{1}_M^T-(p-p^2)I_M.
\end{eqnarray}
Here, we are using the notation $\mathbbm{1}_M$ to denote the $M$-column vector with all its entries equal to one.
Likewise, the cross correlation vector  $r_{d\bar{u},{k}}$ in \eqref{eq:op_missing} is given by
\begin{align}
r_{d\bar{u},{k}}&\triangleq \mathbb{E}\,\boldsymbol{d}_k(i)\bar{\boldsymbol{u}}^*_{k,i}\nonumber\\
&=\mathbb{E}\boldsymbol{d}_k(i)(I-\boldsymbol{F}_{k,i}) \boldsymbol{u}^*_{k,i}+\underbrace{\mathbb{E}\boldsymbol{d}_k(i){\boldsymbol{F}_{k,i} \boldsymbol{\xi}^*_{k,i}}}_{=0}\nonumber\\
&=(1-p)r_{du,k}.
\end{align}
We assume that the perturbed matrix $R_{\bar{u},k}$ remains invertible.
Applying the matrix inversion lemma \cite{Sayed_sdaptive_2008} to \eqref{eq:cov_censored},
we can relate $\bar{w}_k^o$ from \eqref{eq:op_missing} to $w_k^o$ from \eqref{eq:opt_est} as follows:
\begin{align}
\label{eq:censorest}
\bar{w}^o_k&
=(I_M-Q_k)(1-p)w^o
\end{align}
in terms of the matrix
\begin{align}
Q_k\triangleq R^{-1}_{u,{k}}R_{r,k}(I_M+R^{-1}_{u,{k}}R_{r,k})^{-1}.
\end{align}
It is observed from \eqref{eq:censorest} that the new estimate is biased relative to $w^o$.
It is also obvious that the mean-square-error that
results from using \eqref{eq:op_missing} is lower bounded by $J_{k,\min}$ from \eqref{eq:min_mse_perfect}, i.e.,
\begin{align}
\label{eq:min_missing_data}
{J}_{\mathrm{miss},k}&\triangleq J_k(\bar{w}_k^o)
\geq J_{k,\min}.
\end{align}

To mitigate the bias problem, we extend a construction used in \cite{abdolee2012diffusion} in a different context and
associate an
alternative (de-)regularized quadratic cost with each agent $k$, defined as follows:
\begin{align}
\label{eq:local_cost3_cen}
J'_k(w)\triangleq \mathbb{E}|\boldsymbol{d}_k(i)-\bar{\boldsymbol{u}}_{k,i}w|^2-\|w\|^2_{ T_k}
\end{align}
where $T_k$ is a symmetric matrix to be chosen. 
The stationary point of \eqref{eq:local_cost3_cen} is easily seen to occur at
the following location, where we again assume that the coefficient matrix
$(R_{\bar{u},k}-T_k)$ is invertible as needed:
 \begin{align}
 \label{eq:opt_est_miss}
\bar{w}^o_k&=(1-p)(R_{\bar{u},k}-T_k)^{-1}r_{du,k}.
\end{align}
Therefore, if we refer to \eqref{eq:cov_censored} and \eqref{eq:R_2}, we observe that if we select
$T_k$ as
\begin{align}
\label{eq:reg_matrix}
T_k&=pR_{u,k}+R_{r,k}\nonumber\\
&=(p-p^2)I_M\odot R_{u,k}+p\sigma^2_{\xi,k}I_M-(p-p^2)R_{u,k}
\end{align}
then the solution $\bar{w}^o_k$ from \eqref{eq:opt_est_miss} will agree
with the unbiased original estimate $w^o_k$ from \eqref{eq:opt_est}, i.e.,
${\bar{w}_k^o=w_k^o}$.

Now note that since $p-p^2$ is nonnegative for $0\leq p< 1$,
the matrix $T_k$ in \eqref{eq:reg_matrix} is
seen to be the difference of two nonnegative definite matrices.
Therefore, $T_k$ is in general \emph{sign} indefinite.
More importantly, if the de-regularization matrix $T_k$ is selected as in \eqref{eq:reg_matrix},
the cost function in \eqref{eq:local_cost3_cen} becomes \emph{strongly} convex with a unique minimizer.
This is because the Hessian matrix of $J_k'(w)$ relative to $w$ is positive-definite since
\begin{align}
\nabla_w^2 J_k'(w)\;=\;2(R_{\bar{u},k}-T_k)=2(1-p)R_{u,k}>0
\end{align}
for $0\leq p <1$.

\section{Adaptive Distributed Strategy}
\label{sec:algorithm}
In this section, we develop a distributed strategy to recover $w^o$
from missing data by relying on the modified local cost \eqref{eq:local_cost3_cen}.
For the remainder of the paper, we represent the covariance matrix of the regression data in the following form.
\begin{assumption}
\label{assu:diagonal}
The covariance matrix $R_{u,k}$ is diagonal, which is satisfied if the entries of the regression vector $\boldsymbol{u}_{k,i}$ are uncorrelated with each other.~$\quad\qquad \qquad \qquad \qquad \qquad \qquad \qquad \qquad \qquad\qquad\quad\blacksquare$

\end{assumption}
Under Assumption \ref{assu:diagonal}, it holds that
\begin{eqnarray}
R_{r,k}&=&-pR_{u,k}+p\sigma_{\xi,k}^2I_M\\
R_{\bar{u},k}&=&(1-p)R_{u,k}+p\sigma_{\xi,k}^2I_M\label{eq:cov_rel}\\
T_k&=& p\sigma_{\xi,k}^2I_M
\end{eqnarray}

\noindent  where $T_k$ is now nonnegative-definite.  We shall first assume that the parameters $p$ and $\sigma^2_{\xi,k}$ are known. Later, we estimate $\sigma_{\xi,k}^2$ from the data, assuming an estimate for
$p$ is available (a condition that is generally satisfied in practice --- see Sec. \ref{eq:house_hould} where experimental settings are described).

%
%

To develop a distributed algorithm, we let $\mathcal{N}_k$ denote the set of neighbors of agent $k$.
The network then seeks to solve:
\begin{align}
\label{eq:global_cos}
\mathop{\mathrm{\min}}\limits_{w\in\mathbb{R}^M}~ \sum_{k=1}^{N}{J}'_k(w).
\end{align}
Following arguments similar to \cite{sayed_2012_diff_net,sayed_2013_magazine}, we can derive the following
modified Adapt-then-Combine (mATC) diffusion
strategy for the case of missing data:
\begin{equation}
\label{eq:diffusion_atc}
\left\{
\begin{array}{l}
\boldsymbol{e}_k(i)=\boldsymbol{d}_{k}(i)- \bar{\boldsymbol{u}}_{k,i}\boldsymbol{w}_{k,i-1}\\
\boldsymbol{\phi}_{k,i}=(1+\mu_kp\sigma_{\xi,k}^2)\boldsymbol{w}_{k,i-1}+\mu_k \bar{\boldsymbol{u}}^*_{k,i}\boldsymbol{e}_k(i)
\\
\boldsymbol{w}_{k,i}=\sum_{\ell\in \mathcal{N}_k}
a_{\ell k}\boldsymbol{\phi}_{\ell,i}
\end{array}
\right.
\end{equation}
where $\mu_k$ is a small step-size parameter
and the convex combination coefficients $\{a_{\ell k}\}$ satisfy~\cite{sayed_2012_diff_net,sayed_2013_magazine,cattivelli2010diffusion}:
\begin{align}
a_{\ell k}\geq 0,\quad \sum_{\ell\in \mathcal{N}_k}a_{\ell k}=1,\quad a_{\ell k}=0 ~\text{if}~ \ell \notin \mathcal{N}_k.
\end{align}

\subsection{Estimation of Regularization Parameter}
The distributed algorithm \eqref{eq:diffusion_atc} requires knowledge of the censoring noise variance, $\sigma_{\xi,k}^2$.
We now suggest one way to estimate this noise power.
From \eqref{eq:relation_data}, we write
\begin{align}
\label{eq:reg_u}
\boldsymbol{u}_{k,i}=\bar{\boldsymbol{u}}_{k,i}+\boldsymbol{u}_{k,i}\boldsymbol{F}_{k,i}-\boldsymbol{\xi}_{k,i} \boldsymbol{F}_{k,i},
\end{align}
and use this relation to re-write the measurement model \eqref{eq:data_model} in terms of $\bar{\boldsymbol{u}}_{k,i}$ as follows:
\begin{align}
\label{eq:new_model}
\boldsymbol{d}_k(i)=\bar{\boldsymbol{u}}_{k,i}w^o+(\boldsymbol{u}_{k,i}-\boldsymbol{\xi}_{k,i}) \boldsymbol{F}_{k,i}w^o+\boldsymbol{v}_k(i).
\end{align}
It is seen from \eqref{eq:new_model} that
\begin{align}
\label{eq:rel_optimalvalue}
\mathbb{E}|\boldsymbol{d}_k(i)-\bar{\boldsymbol{u}}_{k,i}w^o|^2=\mathbb{E}|(\boldsymbol{u}_{k,i}-\boldsymbol{\xi}_{k,i}) \boldsymbol{F}_{k,i}w^o|^2+\sigma^2_{v,k},
\end{align}
and, hence
\begin{align}
\label{eq:local_cost_cen}
J_{k,\min}&\triangleq \sigma_{v,k}^2\nonumber\\
&=\mathbb{E}|\boldsymbol{d}_k(i)-{\bar{\boldsymbol{u}}_{k,i}}w^o|^2-\mathbb{E}|(\boldsymbol{u}_{k,i}-\boldsymbol{\xi}_{k,i}) \boldsymbol{F}_{k,i}w^o|^2\nonumber\\
&=\mathbb{E}|\boldsymbol{d}_k(i)-{\bar{\boldsymbol{u}}_{k,i}}w^o|^2-p\|w^o\|^2_{R_{u,k}}-p\sigma_{\xi,k}^2\|w^o\|^2.
\end{align}
After a sufficient number of iterations, and for sufficiently small step-sizes,
the estimate $\boldsymbol{w}_{k,i}$ in \eqref{eq:diffusion_atc} is expected to approach
the global minimizer of \eqref{eq:global_cos}, which we already know is
$w^o$~\cite{sayed_2012_diff_net}. If we therefore replace $w^o$
by $\boldsymbol{w}_{k,i-1}$ in \eqref{eq:local_cost_cen} we get for $i\gg 1$:
\begin{align}
\label{eq:cos_opt}
J_{k,\min}\approx
\mathbb{E}|\boldsymbol{e}_k(i)|^2-p\|\boldsymbol{w}_{k,i-1}\|^2_{ R_{u,k}}-p\sigma_{\xi,k}^2\|\boldsymbol{w}_{k,i-1}\|^2.
\end{align}
It is still not possible to estimate $\sigma_{\xi,k}^2$ directly from \eqref{eq:cos_opt} since
the expression depends on $R_{u,k}$ and $p$.
Suppose, as indicated earlier, that an estimate for $p$ is available, say, $\widehat{p}<1$.
This is a reasonable assumption since in many situations in practice, it
is generally known what percentage of the data is corrupted (as illustrated by the
examples in Sec. V). From \eqref{eq:cov_rel}, we can then  estimate $R_{u,k}$ by writing:
\begin{align}
\label{eq:R_u}
R_{u,k}\approx \frac{1}{1-\widehat{p}}R_{\bar{u},k}-\frac{\widehat{p}}{1-\widehat{p}}\sigma_{\xi,k}^2I_M.
\end{align}
Substituting into \eqref{eq:cos_opt} and solving for an estimate for $\sigma_{\xi,k}^2$ we obtain
\begin{eqnarray}
\widehat{\sigma}_{\xi,k}^2
&{\approx}& \frac{(1-\widehat{p})\mathbb{E}|\boldsymbol{e}_k(i)|^2-\widehat{p}\|\boldsymbol{w}_{k,i-1}\|^2_{{R}_{\bar{u},k}}-(1-\widehat{p})\sigma^2_{v,k}}{\widehat{p}(1-2\widehat{p})\|\boldsymbol{w}_{k,i-1}\|^2}\nonumber\\
&\mathop{\approx}\limits^{(a)}& \frac{(1-\widehat{p})\mathbb{E}|\boldsymbol{e}_k(i)|^2-\widehat{p}\|\boldsymbol{w}_{k,i-1}\|^2_{{R}_{\bar{u},k}}}{\widehat{p}(1-2\widehat{p})\|\boldsymbol{w}_{k,i-1}\|^2}
\end{eqnarray}
where in $(a)$ we assumed that the noise variance, $\sigma^2_{v,k}$, is sufficiently small compared to $\|\boldsymbol{w}_{k,i-1}\|^2_{{R}_{\bar{u},k}}$.
Since $\mathbb{E}|\boldsymbol{e}_k(i)|^2$ and the diagonal matrix $R_{\bar{u},k}$ are unknown, we estimate them by means of smoothing filters from data realizations:
\begin{align}
\widehat{\boldsymbol{R}}_{\bar{u},k}(i)&=(1-\alpha_1)\widehat{\boldsymbol{R}}_{\bar{u},k}(i-1)+\alpha_1 (\bar{\boldsymbol{u}}_{k,i}^*\bar{\boldsymbol{u}}_{k,i})\odot I_M \label{eq:cov_up}\\
\widehat{\boldsymbol{\sigma}}_{e,k}(i)&=(1-\alpha_2)\widehat{\boldsymbol{\sigma}}_{e,k}(i-1)+\alpha_2 |\boldsymbol{e}_k(i)|^2\\
\boldsymbol{g}_k(i)&=\frac{(1-\widehat{p})\widehat{\boldsymbol{\sigma}}_{e,k}(i)-\widehat{p}\|
\boldsymbol{w}_{k,i-1}\|^2_{\widehat{\boldsymbol{R}}_{\bar{u},k}(i)}}{\widehat{p}(1-2\widehat{p})\|\boldsymbol{w}_{k,i-1}\|^2}\label{eq:estimat_var}\\
\widehat{\boldsymbol{\sigma}}_{\xi,k}^2(i)&=(1-\alpha_3)\widehat{\boldsymbol{\sigma}}_{\xi,k}^2(i-1)+\alpha_3\boldsymbol{g}_k(i)\label{eq:var_estimate}
\end{align}
where $0< \alpha_m \ll 1,$ for $m=1,2,3$.
To prevent large fluctuations in estimating $\widehat{\boldsymbol{\sigma}}_{\xi,k}^2(i)$, we also use a smoothing filter for updating
$\widehat{\boldsymbol{\sigma}}_{\xi,k}^2(i)$ in \eqref{eq:var_estimate}. Since the covariance matrix of the regressor is assumed to be diagonal, we use the Hadamard product in \eqref{eq:cov_up}.
It is noted that the algorithm does not require knowledge about the statistics of the data, e.g., the correlation matrices are assumed to be unknown.
The resulting diffusion algorithm, henceforth called modified ATC (mATC), is listed in
Algorithm \ref{algo}.\\
%

\begin{algorithm}
\caption{Modified diffusion algorithm (mATC) for missing data}\label{algo}
\begin{algorithmic}[0]
\STATE
\begin{align}
\boldsymbol{e}_k(i)&\triangleq \boldsymbol{d}_k(i)-{\bar{\boldsymbol{u}}_{k,i}}\boldsymbol{w}_{k,i-1}\\
\boldsymbol{\phi}_{k,i}&=(1+\mu_kp\widehat{\boldsymbol{\sigma}}_{\xi,k}^2(i-1))\boldsymbol{w}_{k,i-1}+\mu_k \bar{\boldsymbol{u}}^*_{k,i} \boldsymbol{e}_k(i)\label{eq:eq1}\\
\boldsymbol{w}_{k,i}&=\sum_{\ell\in \mathcal{N}_k}
a_{\ell k}\boldsymbol{\phi}_{\ell,i}\label{eq:eq2}\\
\widehat{\boldsymbol{R}}_{\bar{u},k}(i)&=(1-\alpha_1)\widehat{\boldsymbol{R}}_{\bar{u},k}(i-1)+\alpha_1 (\bar{\boldsymbol{u}}_{k,i}^*\bar{\boldsymbol{u}}_{k,i})\odot I_M \label{eq:var_es_1}\\
\widehat{\boldsymbol{\sigma}}_{e,k}(i)&=(1-\alpha_2)\widehat{\boldsymbol{\sigma}}_{e,k}(i-1)+\alpha_2 |\boldsymbol{e}_k(i)|^2\label{eq:var_es_2}\\
\boldsymbol{g}_k(i)&=\frac{(1-\widehat{p})\widehat{\boldsymbol{\sigma}}_{e,k}(i)-\widehat{p}\|
\boldsymbol{w}_{k,i-1}\|^2_{\widehat{\boldsymbol{R}}_{\bar{u},k}(i)}}{\widehat{p}(1-2\widehat{p})\|\boldsymbol{w}_{k,i-1}\|^2}\label{eq:var_es_3}\\
\widehat{\boldsymbol{\sigma}}_{\xi,k}^2(i)&=(1-\alpha_3)\widehat{\boldsymbol{\sigma}}_{\xi,k}^2(i-1)+\alpha_3\boldsymbol{g}_k(i)\label{eq:var_es_4}
\end{align}
\end{algorithmic}
\end{algorithm}

It is clear from the listing of the algorithm that the operation of the diffusion strategy \eqref{eq:eq1}--\eqref{eq:eq2} is
coupled to steps \eqref{eq:var_es_1}--\eqref{eq:var_es_4}
for estimating $\sigma_{\xi,k}^2$. Proper operation of the algorithm requires a reliable estimate for $\sigma_{\xi,k}^2$.
There are at least two ways to assist in this regard. One way is to use a small step-size $\alpha_3$ in \eqref{eq:var_es_4}.
A second way is to first run a few iterations of the standard diffusion algorithm without bias correction, i.e.,
\begin{eqnarray}
\label{eq:trad_diff}
\boldsymbol{\phi}_{k,i}&=&\boldsymbol{w}_{k,i-1}+\mu_k \bar{\boldsymbol{u}}^*_{k,i} [\boldsymbol{d}_{k}(i)- \bar{\boldsymbol{u}}_{k,i}\boldsymbol{w}_{k,i-1}]
\nonumber\\
\boldsymbol{w}_{k,i}&=&\sum_{\ell\in \mathcal{N}_k}
a_{\ell k}\boldsymbol{\phi}_{\ell,i}
\end{eqnarray}
and then switch to the censored version shown in the above listing.
In the simulations section we illustrate both scenarios.

\begin{remark}
Since the estimate for $\widehat{\boldsymbol{\sigma}}_{\xi,k}^2(i)$ in \eqref{eq:var_es_4} needs to be nonnegative, we can modify \eqref{eq:var_es_3} to
\begin{align}
\boldsymbol{g}_k(i) = \max \left\{ \frac{(1-\widehat{p})\widehat{\boldsymbol{\sigma}}_{e,k}(i)-\widehat{p}\|
\boldsymbol{w}_{k,i-1}\|^2_{\widehat{\boldsymbol{R}}_{\bar{u},k}(i)}}{\widehat{p}(1-2\widehat{p})\|\boldsymbol{w}_{k,i-1}\|^2},  0 \right \}.
\end{align}
\end{remark}

\section{PERFORMANCE ANALYSIS}
\label{eq:err_ana}


%
\subsection{Error Dynamics}
We associate with each agent the error vectors
\begin{align}
\widetilde{\boldsymbol{\phi}}_{k,i}&\triangleq w^o-\boldsymbol{\phi}_{k,i}\\
 \widetilde{\boldsymbol{w}}_{k,i}&\triangleq w^o-\boldsymbol{w}_{k,i}.
\end{align}
Now, if we subtract $w^o$ from both sides of \eqref{eq:eq1} and \eqref{eq:eq2} and replace
$\boldsymbol{d}_{k}(i)$ by \eqref{eq:new_model},
we obtain
\begin{align}
\widetilde{\boldsymbol{\phi}}_{k,i}=&(1+\mu_kp\widehat{\boldsymbol{\sigma}}_{\xi,k}^2(i-1))\widetilde{\boldsymbol{w}}_{k,i-1}-\mu_kp\widehat{\boldsymbol{\sigma}}_{\xi,k}^2(i-1)w^o-\nonumber\\
&\mu_k \bar{\boldsymbol{u}}^*_{k,i} (\bar{\boldsymbol{u}}_{k,i} \widetilde{\boldsymbol{w}}_{k,i-1}+(\boldsymbol{u}_{k,i}-\boldsymbol{\xi}_{k,i}) \boldsymbol{F}_{k,i}w^o+\boldsymbol{v}_k(i))\label{eq:err_rec_1}\\
\widetilde{\boldsymbol{w}}_{k,i}=&\sum_{\ell\in\mathcal{N}_k}a_{\ell k}\widetilde{\boldsymbol{\phi}}_{\ell,i}.\label{eq:err_rec_2}
\end{align}
We collect the errors from across the network into the block vectors:
\begin{align}
\widetilde{\boldsymbol{\phi}}_i&\triangleq
\mathrm{col}\{\widetilde{\boldsymbol{\phi}}_{1,i},\cdots,\widetilde{\boldsymbol{\phi}}_{N,i}\}\label{eq:block_errors_1}\\
\widetilde{\boldsymbol{w}}_i&\triangleq
\mathrm{col}\{\widetilde{\boldsymbol{w}}_{1,i},\cdots,\widetilde{\boldsymbol{w}}_{N,i}
\}\label{eq:block_errors_2}
\end{align}
and note from \eqref{eq:err_rec_1}--\eqref{eq:err_rec_2} that they satisfy the
following recursions:
\begin{align}
\widetilde{\boldsymbol{\phi}}_i\,=\,&[I_{NM}-\mathcal{M}(\boldsymbol{\bar{\mathcal{R}}}_{i}-p\boldsymbol{\mathcal{K}}_{i-1}))]\widetilde{\boldsymbol{w}}_{i-1}-\mathcal{M}\boldsymbol{s}_i-\nonumber\\
&\mathcal{M}(\boldsymbol{\mathcal{R}}_{e,i} +p \boldsymbol{\mathcal{K}}_{i-1})w^o_e\label{eq:err_rec_2_1}\\
\widetilde{\boldsymbol{w}}_i\,=\,&\mathcal{A}^T\widetilde{\boldsymbol{\phi}}_i\label{eq:err_rec_2_2}
\end{align}
where we introduced the quantities:
\begin{align}
w^o_e&\triangleq \mathbbm{1}_N\otimes w^o\label{eq:1}\\
\mathcal{A}&\triangleq A\otimes I_M\label{eq:2}\\
\boldsymbol{\bar{\mathcal{R}}}_{i}&\triangleq \mathrm{diag}\{ \bar{\boldsymbol{u}}^*_{1,i}\bar{\boldsymbol{u}}_{1,i}, \bar{\boldsymbol{u}}^*_{2,i}\bar{\boldsymbol{u}}_{2,i},\ldots, \bar{\boldsymbol{u}}^*_{N,i}\bar{\boldsymbol{u}}_{N,i}\}\label{eq:3} \\
{\boldsymbol{\mathcal{R}}}_{e,i}&\triangleq \mathrm{diag}\left\{ \{\bar{\boldsymbol{u}}^*_{k,i}({\boldsymbol{u}}_{k,i}-\boldsymbol{\xi}_{k,i})\boldsymbol{F}_{k,i}\}_{k=1,\ldots,N}\right\}\label{eq:4}\\
\boldsymbol{\mathcal{K}}_{i-1}&\triangleq\mathrm{diag}\left\{\widehat{\boldsymbol{\sigma}}_{\xi,1}^2(i-1)I_M,\ldots,\widehat{\boldsymbol{\sigma}}_{\xi,N}^2(i-1)I_M\right\}\label{eq:5}\\
\mathcal{M}&\triangleq\mathrm{diag}\{\mu_1I_M,\mu_2I_M,\ldots,\mu_N I_N\}\label{eq:6}\\
\boldsymbol{s}_i&\triangleq \mathrm{col}\{\bar{\boldsymbol{u}}^*_{1,i}\boldsymbol{v}_1(i),\ldots,\bar{\boldsymbol{u}}^*_{N,i}\boldsymbol{v}_N(i)\}\label{eq:7}
\end{align}
where the matrix $A$ is left-stochastic, i.e., $A^T\mathbbm{1}_M=\mathbbm{1}_M$, with its $(\ell,k)$ entry equal to $a_{\ell k}$.
If we now combine \eqref{eq:err_rec_2_1} and \eqref{eq:err_rec_2_2} we find
that $\widetilde{\boldsymbol{w}}_i$ evolves according to the following dynamics:
\begin{align}
\label{eq:err_rec_3}
\widetilde{\boldsymbol{w}}_i&=\mathcal{A}^T[I_{NM}-\mathcal{M}(\boldsymbol{\bar{\mathcal{R}}}_{i}-p\boldsymbol{\mathcal{K}}_{i-1})]\widetilde{\boldsymbol{w}}_{i-1}-\mathcal{A}^T\mathcal{M}\boldsymbol{s}_i\nonumber\\
&\quad-\mathcal{A}^T\mathcal{M}(\boldsymbol{\mathcal{R}}_{e,i} +p \boldsymbol{\mathcal{K}}_{i-1})w^o_e
\end{align}

From the definitions in \eqref{eq:3}, \eqref{eq:4}, and \eqref{eq:7}, we get
\begin{align}
\label{eq:mean_vect}
\mathbb{E} \boldsymbol{s}_i&=0\\
\mathcal{S}&\triangleq \mathbb{E} \boldsymbol{s}_i \boldsymbol{s}_i^*=\mathrm{diag}\{\sigma^2_{v,1}R_{\bar{u},1},\ldots,\sigma^2_{v,N}R_{\bar{u},N}\}\\
\bar{\mathcal{R}}&\triangleq\mathbb{E}\boldsymbol{\bar{\mathcal{R}}}_{i}\nonumber\\
&=\mathrm{diag}\{R_{\bar{u},1},\ldots, R_{\bar{u},N} \}\nonumber\\
&=(1-p)\mathcal{R}+p\,\mathrm{diag}\{\sigma_{\xi,1}^2 I_{M},\ldots,\sigma_{\xi,N}^2 I_{M}\}\label{eq:R_bar}\\
\mathcal{R}_{e}&\triangleq \mathbb{E}\boldsymbol{\mathcal{R}}_{e,i}=-p\,\mathrm{diag}\{\sigma_{\xi,1}^2 I_{M},\ldots,\sigma_{\xi,N}^2 I_{M}\}\label{eq:R_e_bar}
\end{align}
where
\begin{align}
\label{eq:cov_orginial_reg_network}
\mathcal{R}\triangleq \mathrm{diag}\{{R}_{u,1},\ldots,{R}_{u,N}\}
\end{align}
and where we used the following result to compute $\mathcal{R}_{e}$ in \eqref{eq:R_e_bar}:
\begin{align}
&\mathbb{E} \bar{\boldsymbol{u}}^*_{k,i}({\boldsymbol{u}}_{k,i}-\boldsymbol{\xi}_{k,i})\boldsymbol{F}_{k,i}\nonumber\\
&\quad =\mathbb{E}\{
[(I_M-\boldsymbol{F}_{k,i})\boldsymbol{u}_{k,i}^*+\boldsymbol{F}_{k,i}\boldsymbol{\xi}^*_{k,i}]
(\boldsymbol{u}_{k,i}-\boldsymbol{\xi}_{k,i})\boldsymbol{F}_{k,i}\}\nonumber\\
&\quad=\mathbb{E}\{
\underbrace{(I_M-\boldsymbol{F}_{k,i})R_{u,k}\boldsymbol{F}_{k,i}}_{=0}\}
-\sigma^2_{\xi,k}\mathbb{E}\,\boldsymbol{F}^2_{k,i}\nonumber\\
&\quad=-p\sigma^2_{\xi,k}I_M.
\end{align}
\subsection{Long-Term Approximations}
It is clear from \eqref{eq:var_es_1}--\eqref{eq:var_es_4} that the operation of the diffusion strategy \eqref{eq:eq1}--\eqref{eq:eq2} is coupled
with the estimation of the noise power $\sigma_{\xi,k}^2$. This is because
the estimate $\widehat{\boldsymbol{\sigma}}_{\xi,k}^2(i)$ in \eqref{eq:var_es_4} is dependent on $\boldsymbol{w}_{k,i-1}$.
This coupling makes the performance analysis of the algorithm rather challenging.
Since we are mainly interested in assessing the performance of the solution in the infinite-horizon after
sufficient iterations have elapsed, and after the algorithm has been given sufficient
time to learn, we are going to proceed from this point onwards under the assumption
that $i\gg 1$ and that the smoothing filters \eqref{eq:var_es_1}--\eqref{eq:var_es_4} have approached steady-state operation. Specifically, in steady state,
i.e., for $i\rightarrow \infty$, taking the expectation of both sides of Eq.\,\eqref{eq:var_es_1} and considering
$(1-\alpha_1)^i\rightarrow 0$ for $i\gg1$, we obtain:
\begin{align}
\mathbb{E} \widehat{\boldsymbol{R}}_{\bar{u},k}(i)&=\sum_{j=0}^i\alpha_1(1-\alpha_1)^{i-j}\mathbb{E} (\bar{\boldsymbol{u}}_{k,{j}}^*\bar{\boldsymbol{u}}_{k,{j}})\odot I_M\nonumber\\
&=\alpha_1\frac{1-(1-\alpha_1)^{i+1}}{1-(1-\alpha_1)}R_{\bar{u},k}\odot I_M\nonumber\\
& \rightarrow R_{\bar{u},k},\quad i\gg 1
\end{align}
so that $\widehat{\boldsymbol{R}}_{\bar{u},k}$ tends on average to the true value $R_{\bar{u},k}$.
Similarly, for  Eqs. \eqref{eq:var_es_2}--\eqref{eq:var_es_4}, they approach steady-state operation with
\begin{align}
\mathbb{E}\widehat{\boldsymbol{\sigma}}_{e,k}^2(i)&\rightarrow \mathbb{E}|\boldsymbol{e}_k(i)|^2, \quad i\gg 1\label{eq:st_st_1}\\
\mathbb{E} \widehat{\boldsymbol{\sigma}}_{\xi,k}^2(i)&\rightarrow \mathbb{E}\boldsymbol{g}_k(i),\quad i\gg 1.\label{eq:st_st_2}
\end{align}
We now estimate $ \mathbb{E} \boldsymbol{g}_k(i)$ by employing a first-order Taylor series approximation.
Let $\boldsymbol{z}=[\boldsymbol{z}_1,\ldots,\boldsymbol{z}_n]$ be a real random vector with
mean $\mathbb{E} \boldsymbol{z}=[\mathbb{E}\boldsymbol{z}_1,\ldots,\mathbb{E}\boldsymbol{z}_n]$.
The first-order Taylor series expansion of a differentiable real function $f(\boldsymbol{z})$ about the mean $\mathbb{E} \boldsymbol{z}$
can be expressed as follows \cite[p.~241]{2002statistical}--\cite[p.~295]{Kay_93}:
\begin{align}
\label{eq:taylor_series}
f(\boldsymbol{z})\approx f(\mathbb{E}\boldsymbol{z})+
\sum_{k=1}^n\frac{\partial}{\partial \boldsymbol{z}_k}f(\boldsymbol{z}_k)|_{\boldsymbol{z}_k=\mathbb{E} \boldsymbol{z}_k}(\boldsymbol{z}_k-\mathbb{E} \boldsymbol{z}_k).
\end{align}
If we evaluate the expectation of both sides of \eqref{eq:taylor_series} with respect to the random vector $\boldsymbol{z}$,
we get the follwoing approximation:
 \begin{align}
\label{eq:taylor_series_expect}
\mathbb{E}f(\boldsymbol{z})\approx f(\mathbb{E}\boldsymbol{z}).
\end{align}
Now consider a function of the ratio of two random variables as $f(\boldsymbol{x},\boldsymbol{y})=\boldsymbol{x}/\boldsymbol{y}$, and assume that $\boldsymbol{y}$
has nonzero mean. From \eqref{eq:taylor_series_expect}, we can write
\begin{align}
\label{eq:expec_ratio}
\mathbb{E}f(\boldsymbol{x},\boldsymbol{y})\approx\frac{\mathbb{E}\boldsymbol{x}}{\mathbb{E}\boldsymbol{y}}.
\end{align}
To approximate $ \mathbb{E} \boldsymbol{g}_k(i)$, we apply \eqref{eq:expec_ratio} to \eqref{eq:var_es_3} as follows:
\begin{align}
\mathbb{E} \boldsymbol{g}_k(i)\approx \frac{(1-\widehat{p})\mathbb{E}\widehat{\boldsymbol{\sigma}}_{k,e}^2(i)-\widehat{p}\,\mathbb{E}
\left(\|\boldsymbol{w}_{k,i}\|^2_{\widehat{\boldsymbol{R}}_{\bar{u},k}(i)}\right)}{\widehat{p}(1-2\widehat{p})\mathbb{E}\|\boldsymbol{w}_{k,i}\|^2}.
\end{align}

For tractability, we assume that, in steady-state:
\begin{align}
\mathbb{E}\left(\|
\boldsymbol{w}_{k,i}\|^2_{\widehat{\boldsymbol{R}}_{\bar{u},k}(i)}\right)\approx  \|
{w}^o\|^2_{R_{\bar{u},k}}.
\end{align}
so that
\begin{align}
\mathbb{E} \boldsymbol{g}_k(i)&\approx \frac{(1-\widehat{p})\mathbb{E} |\boldsymbol{e}_k(i)|^2-\widehat{p}\|
w^o\|^2_{R_{\bar{u},k}}}{\widehat{p}(1-2\widehat{p})\|w^o\|^2}\nonumber\\
&\mathop{\approx}\limits^{(a)} \frac{(1-\widehat{p})(\mathbb{E} |\boldsymbol{e}_k(i)|^2-\widehat{p}\|
w^o\|^2_{R_{u,k}})-\widehat{p}^2\sigma^2_{\xi,k}\|w^o\|^2}{\widehat{p}(1-2\widehat{p})\|w^o\|^2}\nonumber\\
&\mathop{\approx}\limits^{(b)} \frac{(1-\widehat{p})(\sigma^2_{v,k}+\widehat{p}\sigma^2_{\xi,k}\|w^o\|^2)-\hat{p}^2\sigma^2_{\xi,k}\|w^o\|^2}{\widehat{p}(1-2\widehat{p})\|w^o\|^2}\nonumber\\
&\mathop{\approx}\limits^{(c)}\sigma^2_{\xi,k}
\end{align}
where in $(a)$ we replace $R_{\bar{u},k}$ from \eqref{eq:cov_rel}, in $(b)$ we use the relation from \eqref{eq:local_cost_cen}, and in $(c)$ we
assume that the term $(1-\widehat{p})\sigma^2_{v,k}$ is sufficiently small compared to $\widehat{p}(1-2\widehat{p})\sigma^2_{\xi,k}\|w^o\|^2$.
Therefore, in steady state, we set
\begin{align}
\label{eq:mean_var_estimate}
\mathcal{K}\triangleq \lim_{i\rightarrow \infty}\mathbb{E}
\,\boldsymbol{\mathcal{K}}_i\approx \mathrm{diag}\left\{\sigma_{\xi,1}^2I_M,\ldots,\sigma_{\xi,N}^2I_M\right\}.
\end{align}

%

\subsection{Mean Stability Analysis}
First note that, in steady state, from \eqref{eq:R_e_bar} and \eqref{eq:mean_var_estimate} we have
\begin{align}
\label{eq:tmp}
\lim_{i\rightarrow \infty}\mathbb{E}(\boldsymbol{\mathcal{R}}_{e,i} +p \boldsymbol{\mathcal{K}}_{i-1})\approx 0.
\end{align}
Now, since the variables $\boldsymbol{u}_{k,i}$ and $\boldsymbol{\xi}_{k,i}$ are temporally  white and spatially independent, then the error vectors $\widetilde{\boldsymbol{w}}_{\ell,j}$ are independent of $\boldsymbol{u}_{k,i}$ and $\boldsymbol{\xi}_{k,i}$ for all $j$ if $k\neq \ell$ and for $k=\ell$ if $j\leq i-1$.
Therefore, taking expectation of both sides of \eqref{eq:err_rec_3} for large enough $i$ gives under the long-term approximations of the previous section:
\begin{align}
\label{eq:recur}
\mathbb{E}\widetilde{\boldsymbol{w}}_i=\mathcal{A}^T[I_{NM}-(1-p)\mathcal{M}\mathcal{R}]\mathbb{E}\widetilde{\boldsymbol{w}}_{i-1},\quad i\gg1.
\end{align}
This recursion is stable if the step sizes are sufficiently small and satisfy
\begin{align}
\label{eq:stepsize_cons}
0<\mu_k <
\frac{2}{(1-p)\lambda_{\max}(R_{u,k})}
\end{align}
where $\lambda_{\max}(\cdot)$ denotes the maximum eigenvalue of its matrix argument.
In this case, the estimator becomes asymptotically unbiased, i.e., $\lim_{i\rightarrow \infty}\mathbb{E}\widetilde{\boldsymbol{w}}_i= 0$.

\subsection{Mean Variance Analysis}
We rewrite \eqref{eq:err_rec_3} more compactly as
\begin{align}
\label{eq:err_rec_4}
\widetilde{\boldsymbol{w}}_i=\boldsymbol{\mathcal{B}}_i\widetilde{\boldsymbol{w}}_{i-1}-\mathcal{\mathcal{G}}\boldsymbol{s}_i-\boldsymbol{\mathcal{D}}_iw^o_e
\end{align}
where
\begin{align}
\boldsymbol{\mathcal{B}}_i&\triangleq \mathcal{A}^T[I_{NM}-\mathcal{M}({\boldsymbol{\bar{\mathcal{R}}}}_{i}-p\boldsymbol{\mathcal{K}}_{i-1})]\label{eq:B_ex}\\
\boldsymbol{\mathcal{D}}_i&\triangleq \mathcal{A}^T\mathcal{M}[\boldsymbol{\mathcal{R}}_{e,i}+p\boldsymbol{\mathcal{K}}_{i-1}]\label{eq:D_ex}\\
\mathcal{G}&\triangleq \mathcal{A}^T\mathcal{M}\label{eq:def_G}.
\end{align}
\noindent The mean-square error analysis of the algorithm relies on evaluating a weighted variance of the error vector. Let $\Sigma$ denote an
arbitrary nonnegative-definite matrix that we are free to choose. We express the weighted square measures on both sides of \eqref{eq:err_rec_4} as
\begin{align}
\label{eq:norm_rel}
\|\widetilde{\boldsymbol{w}}_i\|^2_{\Sigma}=&\|\boldsymbol{\mathcal{B}}_i\widetilde{\boldsymbol{w}}_{i-1}-\mathcal{G}\boldsymbol{s}_i-\boldsymbol{\mathcal{D}}_iw^o_e\|^2_{\Sigma}\nonumber\\
=&\widetilde{\boldsymbol{w}}^*_{i-1}\boldsymbol{\mathcal{B}}^*_i\Sigma\boldsymbol{\mathcal{B}}_i\widetilde{\boldsymbol{w}}_{i-1}+\boldsymbol{s}^*_i\mathcal{G}^T\Sigma\mathcal{G}\boldsymbol{s}_i+{w_e^{o*}}\boldsymbol{\mathcal{D}}^*_i\Sigma\boldsymbol{\mathcal{D}}_iw^o_e\nonumber\\
&-\widetilde{\boldsymbol{w}}_{i-1}^*\boldsymbol{\mathcal{B}}_i^*\Sigma\mathcal{G}\boldsymbol{s}_i-\boldsymbol{s}^*_i\mathcal{G}^T\Sigma\boldsymbol{\mathcal{B}}_i\widetilde{\boldsymbol{w}}_{i-1}\nonumber\\
&-\widetilde{\boldsymbol{w}}_{i-1}^*\boldsymbol{\mathcal{B}}_i^*\Sigma\boldsymbol{\mathcal{D}}_iw^o_e-w_e^{o*}\boldsymbol{\mathcal{D}}^*_i\Sigma\boldsymbol{\mathcal{B}}_i\widetilde{\boldsymbol{w}}_{i-1}\nonumber\\
&+{\boldsymbol{s}}_{i}^*\mathcal{G}\Sigma\boldsymbol{\mathcal{D}}_iw^o_e+w_e^{o*}\boldsymbol{\mathcal{D}}^*_i\Sigma\mathcal{G}{\boldsymbol{s}}_{i}.
\end{align}
We now compute the expectation of both sides of \eqref{eq:norm_rel} in steady state.
Since $\lim_{i\rightarrow \infty}\mathbb{E}\widetilde{\boldsymbol{w}}_i= 0$,
$\mathbb{E} \boldsymbol{s}_i=0$, and $\widetilde{\boldsymbol{w}}_{i-1}$ and $\boldsymbol{s}_i$ are independent of each other, we get the following results:
\begin{align}
\mathbb{E}\widetilde{\boldsymbol{w}}^*_{i-1}\boldsymbol{\mathcal{B}}^*_i\Sigma\mathcal{G}\boldsymbol{s}_i&=0\\
\mathbb{E}{w^o_e}^*\boldsymbol{\mathcal{D}}^*_i\Sigma\mathcal{G}\boldsymbol{s}_i&=0\\
\lim_{i\rightarrow \infty}\mathbb{E}\widetilde{\boldsymbol{w}}^*_{i-1}\boldsymbol{\mathcal{B}}_i\Sigma\boldsymbol{\mathcal{D}}_iw^o_e
&=\lim_{i\rightarrow \infty}\mathbb{E}[\mathbb{E}(\widetilde{\boldsymbol{w}}^*_{i-1}\boldsymbol{\mathcal{B}}^*_i\Sigma\boldsymbol{\mathcal{D}}_i|_{\widetilde{\boldsymbol{w}}_{i-1}})]w^o_e\nonumber\\
&=\lim_{i\rightarrow \infty}\mathbb{E}\widetilde{\boldsymbol{w}}^*_{i-1}[\mathbb{E}(\boldsymbol{\mathcal{B}}^*_i\Sigma\boldsymbol{\mathcal{D}}_i|_{\widetilde{\boldsymbol{w}}_{i-1}})]w^o_e\nonumber\\
&=\lim_{i\rightarrow \infty}\mathbb{E}\widetilde{\boldsymbol{w}}^*_{i-1}[\mathbb{E}\mathcal{A}\Sigma\boldsymbol{\mathcal{D}}_i]w^o_e+O(\mathcal{M}^2)\nonumber\\
&\approx\lim_{i\rightarrow \infty}\mathbb{E}\widetilde{\boldsymbol{w}}^*_{i-1}[\mathbb{E}\mathcal{A}\Sigma\boldsymbol{\mathcal{D}}_i]w^o_e\nonumber\\
&\approx 0, \label{eq:cross_3}
\end{align}
since from \eqref{eq:tmp}
\begin{align}
\mathbb{E} \boldsymbol{\mathcal{D}}_i\approx 0,\quad i\gg 1.
\end{align}
Noting that the cross-terms are either zero or negligible under expectation, we conclude from \eqref{eq:err_rec_4}--\eqref{eq:cross_3} that
\begin{align}
\label{eq:Ex_norm}
\lim_{i\rightarrow \infty}
\mathbb{E}\|\widetilde{\boldsymbol{w}}_i\|^2_{\Sigma}=&
\lim_{i\rightarrow \infty}
\big[\mathbb{E}(\widetilde{\boldsymbol{w}}^*_{i-1}\boldsymbol{\mathcal{B}}^*_i\Sigma\boldsymbol{\mathcal{B}}_i\widetilde{\boldsymbol{w}}_{i-1})+\mathbb{E}(\boldsymbol{s}^*_i\mathcal{G}^T\Sigma\mathcal{G}\boldsymbol{s}_i)\nonumber\\
&+\mathbb{E}({w_e^{o*}}\boldsymbol{\mathcal{D}}^*_i\Sigma\boldsymbol{\mathcal{D}}_iw^o_e)\big].
\end{align}
We now evaluate the terms that appear on the right-hand side of \eqref{eq:Ex_norm}.
The first term can be written as
\begin{align}
&\lim_{i\rightarrow \infty}
\mathbb{E}(\widetilde{\boldsymbol{w}}^*_{i-1}\boldsymbol{B}^*_i\Sigma\boldsymbol{B}_i\widetilde{\boldsymbol{w}}_{i-1})=\lim_{i\rightarrow \infty}
\mathbb{E}\|\widetilde{\boldsymbol{w}}_{i-1}\|^2_{\Sigma'}
\end{align}
where
\begin{align}
\label{eq:sigma_1}
\Sigma'\triangleq&\lim_{i\rightarrow \infty}
\mathbb{E}(\boldsymbol{\mathcal{B}}^*_i\Sigma\boldsymbol{\mathcal{B}}_i)\nonumber\\
=&\mathcal{A}\Sigma\mathcal{A}^T-(1-p)\mathcal{R}\mathcal{M}\mathcal{A}\Sigma\mathcal{A}^T-(1-p)\mathcal{A}{\Sigma}\mathcal{A}^T\mathcal{M}{\mathcal{R}}\nonumber\\
&+O(\mathcal{M}^2),
\end{align}
in view of the fact that
\begin{align}
\lim_{i\rightarrow \infty} \mathbb{E}(\bar{\boldsymbol{\mathcal{R}}}_i-p \boldsymbol{\mathcal{K}}_{i-1})=(1-p)\mathcal{R}
\end{align}
with $\mathcal{R}$ defined by \eqref{eq:cov_orginial_reg_network}.
For sufficiently small step-sizes, the term $O(\mathcal{M}^2)$ in \eqref{eq:sigma_1} is negligible
and we set
\begin{align}
\label{eq:sig_prim}
\Sigma'=\mathcal{A}\Sigma\mathcal{A}^T-(1-p){\mathcal{R}}\mathcal{M}\mathcal{A}\Sigma\mathcal{A}^T-(1-p)\mathcal{A}{\Sigma}\mathcal{A}^T\mathcal{M}{\mathcal{R}}.
\end{align}
The second term on the right-hand side of \eqref{eq:Ex_norm} can be calculated as
\begin{align}
\label{eq:S-expression}
\mathbb{E}(\boldsymbol{s}^*_i\mathcal{G}^T\Sigma\mathcal{G}\boldsymbol{s}_i)
=\mathbb{E}\, \mathrm{Tr}(\mathcal{G}^T\Sigma\mathcal{G}\boldsymbol{s}_i\boldsymbol{s}_i^*)=
\mathrm{Tr}(\mathcal{G}^T\Sigma\mathcal{G}\mathcal{S}).
\end{align}
Now it is clear from expression \eqref{eq:4} that $\boldsymbol{\mathcal{R}}_{e,i}$ depends on random variables that are available at time $i$, while
expression \eqref{eq:5} shows that $\boldsymbol{\mathcal{K}}_{i-1}$ depends on different random variables up to time $i-1$.
Therefore, $\boldsymbol{\mathcal{R}}_{e,i}$ and $\boldsymbol{\mathcal{K}}_{i-1}$ are independent of each other.
Then, the last term on the right-hand side of \eqref{eq:Ex_norm} can be computed as follows:
\begin{align}
\label{eq:last_term}
&\lim_{i\rightarrow \infty}\mathbb{E}({w_e^{o*}}\boldsymbol{\mathcal{D}}^*_i\Sigma\boldsymbol{\mathcal{D}}_iw^o_e)\nonumber\\
&\qquad=\lim_{i\rightarrow \infty}\big[
2pw_e^{o*}\mathrm{Re}\{\mathbb{E}\boldsymbol{\mathcal{R}}^*_{e,i}\mathcal{MA}\Sigma\mathcal{A}^T\mathcal{M}\boldsymbol{\mathcal{K}}^*_{i-1}\}w_e^{o}\nonumber\\
&\qquad\quad+w_e^{o*}\mathbb{E}(\boldsymbol{\mathcal{R}}^*_{e,i}\mathcal{M}\mathcal{A}\Sigma\mathcal{A}^T\mathcal{M} \boldsymbol{\mathcal{R}}_{e,i} )w_e^{o}\nonumber\\
&\qquad\quad+p^2w_e^{o*}\mathbb{E}(\boldsymbol{\mathcal{K}}^*_{i-1}\mathcal{M}\mathcal{A}\Sigma\mathcal{A}^T\mathcal{M} \boldsymbol{\mathcal{K}}_{i-1})w_e^{o}\big]\nonumber\\
&\qquad\mathop{=}\limits^{(a)}-2w_e^{o*}\mathcal{R}_{e}\mathcal{MA}\Sigma\mathcal{A}^T\mathcal{M}\mathcal{R}_{e}w_e^{o}\nonumber\\
&\qquad\quad+ \lim_{i\rightarrow \infty} w_e^{o*}\mathbb{E}(\boldsymbol{\mathcal{R}}^*_{e,i}\mathcal{M}\mathcal{A}\Sigma\mathcal{A}^T\mathcal{M} \boldsymbol{\mathcal{R}}_{e,i} )w_e^{o}\nonumber\\
&\qquad\quad+\lim_{i\rightarrow \infty}p^2w_e^{o*}\mathbb{E}[\boldsymbol{\mathcal{K}}^*_{i-1}\mathcal{M}\mathcal{A}\Sigma\mathcal{A}^T\mathcal{M} \boldsymbol{\mathcal{K}}_{i-1}]w_e^{o}.
\end{align}
In $(a)$, we used the independence of $\boldsymbol{\mathcal{R}}_{e,i}$ and $\boldsymbol{\mathcal{K}}_{i-1}$ and the
fact from \eqref{eq:tmp} that
$\mathbb{E}\boldsymbol{\mathcal{R}}_{e,i}=-p\mathbb{E}\boldsymbol{\mathcal{K}}_{i-1},~i\gg1$.
In order to obtain a more compact representation for the variance relation, we use the following vector notation:
\begin{align}
\sigma=\mathrm{vec}(\Sigma),\quad \sigma'=\mathrm{vec}(\Sigma'),
\end{align}
where the $\mathrm{vec}$ operator vectorizes a matrix by placing its columns on top of each other.
We also use the following Kroneceker product properties~\cite{sayed_2012_diff_net}:
\begin{align}
&\mathrm{vec}(ABC)=(C^T\otimes A) \mathrm{vec}(B)\nonumber\\
&\mathrm{Tr}(AB)=[\mathrm{vec}(B^T)]^T\mathrm{vec}(A).
\end{align}
Then, from \eqref{eq:sig_prim} we can write
\begin{align}
\sigma'=\mathcal{F}\sigma,
\end{align}
where
\begin{align}
\mathcal{F}\triangleq \mathcal{A}\otimes \mathcal{A}-(1-p)\mathcal{A}\otimes {\mathcal{R}}\mathcal{MA}-(1-p){\mathcal{R}^T}\mathcal{M}\mathcal{A}\otimes \mathcal{A}.
\end{align}
Likewise, we have
\begin{align}
\mathrm{Tr}(\mathcal{G}^T\Sigma\mathcal{G}\mathcal{S})=\mathrm{Tr}(\Sigma\mathcal{G}\mathcal{S}\mathcal{G}^T)=[\mathrm{vec}(\mathcal{GS}^T\mathcal{G}^T)]^T\sigma.
\end{align}
In a similar manner, we can express the right-hand side of \eqref{eq:last_term} using the $\mathrm{vec}$ operator as follows:
\begin{align}
\label{eq:z_term}
&\lim_{i\rightarrow \infty}\mathbb{E}({w_e^{o*}}\boldsymbol{\mathcal{D}}^*_i\Sigma\boldsymbol{\mathcal{D}}_iw^o_e)\nonumber\\
&\quad=-2\mathrm{Tr}(\Sigma\mathcal{A}^T\mathcal{M}\mathcal{R}_e w^o_ew^{o*}_e\mathcal{R}_e\mathcal{MA})\nonumber\\
&\quad\quad+\mathrm{Tr}(\Sigma\mathcal{A}^T\mathcal{M}\lim_{i\rightarrow \infty}\mathbb{E}(\boldsymbol{\mathcal{R}}_{e,i}w^o_ew^{o*}_e\boldsymbol{\mathcal{R}}^*_{e,i})\mathcal{MA})\nonumber\\
&\quad\quad+p^2\mathrm{Tr}(\Sigma\mathcal{A}^T\mathcal{M}\lim_{i\rightarrow \infty}{\mathbb{E}(\boldsymbol{\mathcal{K}}_{i-1}w^o_ew^{o*}_e\boldsymbol{\mathcal{K}}^T_{i-1})}\mathcal{MA})\nonumber\\
&\quad=-2[\mathrm{vec}(\mathcal{A}^T\mathcal{M}\mathcal{R}_e w^o_ew^{o*}_e\mathcal{R}_e\mathcal{MA})^T]^T\sigma\nonumber\\
&\quad\quad+[\mathrm{vec}(\mathcal{A}^T\mathcal{M}\lim_{i\rightarrow \infty}\mathbb{E}(\boldsymbol{\mathcal{R}}_{e,i}w^o_ew^{o*}_e\boldsymbol{\mathcal{R}}^*_{e,i})\mathcal{MA})^T]^T\sigma\nonumber\\
&\quad\quad+p^2[\mathrm{vec}(\mathcal{A}^T\mathcal{M}\lim_{i\rightarrow \infty}\underbrace{\mathbb{E}(\boldsymbol{\mathcal{K}}_{i-1}w^o_ew^{o*}_e\boldsymbol{\mathcal{K}}^T_{i-1})}_{\triangleq \Pi}\mathcal{MA})^T]^T\sigma.
\end{align}

\noindent It is noted that the quantity $\mathbb{E}(\boldsymbol{\mathcal{R}}_{e,i}w^o_ew^{o*}_e\boldsymbol{\mathcal{R}}^*_{e,i})$, in general, does not have
a closed form expression.
We can approximate it by means of ensemble averaging.

The last term on the right-hand side of \eqref{eq:z_term} can be approximated as follows.
Let ${\Omega}\triangleq w^o_ew^{o*}_e$ be a Hermitian block matrix, whose $(n,m)$ block is given by
\begin{align}
\Omega_{n,m}=w^ow^{o*}.
\end{align}
The $(n,m)$-th block of $\Pi$ can be
obtained as follows:
\begin{align}
\Pi_{n,m}&=\mathbb{E}(\widehat{\boldsymbol{\sigma}}^2_{\xi,n}(i-1)\Omega_{n,m}\widehat{\boldsymbol{\sigma}}^2_{\xi,m}(i-1))\nonumber\\
&=\Omega_{n,m}\mathbb{E}(\widehat{\boldsymbol{\sigma}}^2_{\xi,n}(i-1)\widehat{\boldsymbol{\sigma}}^2_{\xi,m}(i-1)).
\end{align}
From \eqref{eq:var_es_4}, we can write
\begin{align}
\label{eq:cross_up}
\widehat{\boldsymbol{\sigma}}_{\xi,n}^2(i)\widehat{\boldsymbol{\sigma}}_{\xi,m}^2(i)=&
(1-\alpha_3)^2\widehat{\boldsymbol{\sigma}}_{\xi,n}^2(i-1)\widehat{\boldsymbol{\sigma}}_{\xi,m}^2(i-1)\nonumber\\
&+\alpha^2_3\boldsymbol{g}_n(i)\boldsymbol{g}_m(i)\nonumber\\
&+\alpha_3(1-\alpha_3)[\widehat{\boldsymbol{\sigma}}_{\xi,n}^2(i-1)\boldsymbol{g}_m(i)\nonumber\\
&+\widehat{\boldsymbol{\sigma}}_{\xi,m}^2(i-1)\boldsymbol{g}_n(i)].
\end{align}
In steady state, we have $\mathbb{E} \widehat{\boldsymbol{\sigma}}_{\xi,n}^2(i)\widehat{\boldsymbol{\sigma}}_{\xi,m}^2(i)
=\mathbb{E} \widehat{\boldsymbol{\sigma}}_{\xi,n}^2(i-1)\widehat{\boldsymbol{\sigma}}_{\xi,m}^2(i-1)$.
We further benefit from the following approximations for $i\gg 1$:
\begin{align}
\mathbb{E}|\boldsymbol{e}_n(j)|^2|\boldsymbol{e}_m(i)|^2&\approx \mathbb{E}|\boldsymbol{e}_n(j)|^2\mathbb{E}|\boldsymbol{e}_m(i)|^2\nonumber\\
\mathbb{E} \|\boldsymbol{w}_{n,i-1}\|^2\|\boldsymbol{w}_{m,i-1}\|^2&\approx \|w^o\|^4\nonumber\\
\mathbb{E}\|\boldsymbol{w}_{n,i-1}\|^2_{{\widehat{\boldsymbol{R}}_{\bar{u},n}}(i)}\|\boldsymbol{w}_{m,i-1}\|^2_{{\widehat{\boldsymbol{R}}_{\bar{u},m}}(i)}&\approx
\|w^o\|^2_{R_{\bar{u},n}}\|w^o\|^2_{R_{\bar{u},m}}.
\end{align}
We again use a first-order Taylor series expansion to approximate the following terms as:
\begin{align}
\mathbb{E}\boldsymbol{g}_n(i)\boldsymbol{g}_m(i)&\approx
\mathbb{E}\boldsymbol{g}_n(i)\mathbb{E}\boldsymbol{g}_m(i)\approx\sigma^2_{\xi,n}\sigma^2_{\xi,m}\\
\mathbb{E}\widehat{\boldsymbol{\sigma}}_{\xi,n}^2(i-1)\boldsymbol{g}_m(i)&\approx
\mathbb{E}\widehat{\boldsymbol{\sigma}}_{\xi,n}^2(i-1)\mathbb{E}\boldsymbol{g}_m(i)\approx\sigma^2_{\xi,n}\sigma^2_{\xi,m}.
\end{align}
Using the above expressions and
\eqref{eq:cross_up}, we obtain
\begin{align}
\label{eq:cross_2}
(2\alpha_3-\alpha_3^2)\mathbb{E}\widehat{\boldsymbol{\sigma}}_{\xi,n}^2(i)\widehat{\boldsymbol{\sigma}}_{\xi,m}^2(i)\approx&
(2\alpha_3-\alpha_3^2)\sigma_{\xi,n}^2\sigma_{\xi,m}^2.
\end{align}


\noindent Hence, in steady state, the matrix $\Pi$ is approximated as
\begin{align}
\Pi\approx \mathcal{K}w^{o}w^{o*}\mathcal{K}
\end{align}
and in steady state, expression \eqref{eq:last_term} is approximated as
\begin{align}
\lim_{i\rightarrow \infty}\mathbb{E}({w_e^{o*}}\boldsymbol{\mathcal{D}}^*_i\Sigma\boldsymbol{\mathcal{D}}_iw^o_e)\approx[\mathrm{vec}(\mathcal{Z}^T)]^T\sigma
\end{align}
where
\begin{align}
\label{eq:Z}
\mathcal{Z}\triangleq&- \mathcal{A}^T\mathcal{M}\mathcal{R}_e w^o_ew^{o*}_e\mathcal{R}_e\mathcal{MA}\nonumber\\
&+\mathcal{A}^T\mathcal{M}\mathbb{E}(\boldsymbol{R}_{e,i}w^o_ew^{o*}_e\boldsymbol{R}^*_{e,i})\mathcal{MA}.
\end{align}
\noindent
Referring back to \eqref{eq:Ex_norm}, and using the notation $\|x\|^2_{\sigma}$
interchangeably with $\|x\|^2_{\Sigma}$, we get
\begin{align}
\label{eq:final_var_relation}
\lim_{i\rightarrow \infty}\mathbb{E}\|\widetilde{\boldsymbol{w}}_i\|_{\sigma}=\lim_{i\rightarrow \infty}\mathbb{E}\|\widetilde{\boldsymbol{w}}_{i-1}\|_{\mathcal{F}\sigma}+[\mathrm{vec}(\mathcal{Z}^T+\mathcal{Y}^T)]^T\sigma
\end{align}
where
\begin{align}
\label{eq:Y_exp}
\mathcal{Y}\triangleq \mathcal{GS}\mathcal{G}^T.
\end{align}
It follows that
\begin{align}
\label{eq:final_var_relation}
\lim_{i\rightarrow \infty}\mathbb{E}\|\widetilde{\boldsymbol{w}}_i\|_{(I-\mathcal{F})\sigma}=[\mathrm{vec}(\mathcal{Z}^T+\mathcal{Y}^T)]^T\sigma.
\end{align}
Based on the variance relation \eqref{eq:final_var_relation}, different quantities can be computed.
For example, we can evaluate the network and individual mean-square deviation (MSDs), respectively, defined as
\begin{align}
\mathrm{MSD}^{\mathrm{network}}&\triangleq \lim_{i\rightarrow \infty}\frac{1}{N}
\sum_{k=1}^N\mathbb{E}\|\widetilde{\boldsymbol{w}}_{k,i}\|^2= \lim_{i\rightarrow \infty}
\mathbb{E}\|\widetilde{\boldsymbol{w}}_{i}\|^2_{\frac{1}{N}}\nonumber\\
\mathrm{MSD}_k&\triangleq  \lim_{i\rightarrow \infty}
\mathbb{E}\|\widetilde{\boldsymbol{w}}_{k,i}\|^2= \lim_{i\rightarrow \infty}
\mathbb{E} \|\widetilde{\boldsymbol{w}}_{i}\|^2_{\mathcal{I}_k}
\end{align}
where $\mathcal{I}_k\triangleq \mathrm{diag}\{0,\ldots,0,I_M,0,\ldots,0\}$, with the identity matrix appearing in the $k-$th block location.

In order to derive the network MSD from variance relation \eqref{eq:final_var_relation}, we select the weighting vector $\sigma$ such that
\begin{align}
(I-\mathcal{F})\sigma=\frac{1}{N}\mathrm{vec}(I_{NM}).
\end{align}
Then, the network MSD can be calculated as
\begin{align}
\mathrm{MSD}^{\mathrm{network}}=\frac{1}{N}[\mathrm{vec}(\mathcal{Z}^T+\mathcal{Y}^T)]^T(I-\mathcal{F})^{-1}\mathrm{vec}(I_{NM}).\label{eq:msd_net}
\end{align}
Likewise, the individual MSD can be computed as
\begin{align}
\mathrm{MSD}_k=[\mathrm{vec}(\mathcal{Z}^T+\mathcal{Y}^T)]^T(I-\mathcal{F})^{-1}\mathrm{vec}(\mathcal{I}_{k}).\label{msd:msd_single}
\end{align}

\section{Design Examples}
\label{sec:Sim_results}
\begin{figure}
\center
\begin{tikzpicture}[scale=2,rotate=-30]
\tikzstyle{every node}=[draw,black,ultra thick,shape=circle,fill=blue!20];
\path (-0.1,-.10) node (v1) {\tiny{1}};
\path (-.2,1.8) node (v2) {\tiny{2}};
\path (.8,.4) node (v3) {\tiny{3}};
\path (1,1.1) node (v4) {\tiny{4}};
\path (1,2) node (v5) {\tiny{5}};
\path (1.8,1) node (v6) {\tiny{6}};
\path (-.7,.9) node (v7) {\tiny{7}};
\draw[black] (v1) -- (v3)
(v1) -- (v7)
(v1) -- (v2)
(v1) -- (v5);
\draw[black] (v2) -- (v3)
(v2) -- (v4)
(v2) -- (v5)
(v2) -- (v7);
\draw[black]
(v3) -- (v4)
(v3) -- (v6)
(v3) -- (v7);
\draw[black]
(v4) -- (v5)
(v4) -- (v6);
\draw[black]
(v5) -- (v6)
(v5) -- (v7);
\end{tikzpicture}
 \caption{Topology of the network used in the simulations.}
\label{fig:network}%
\end{figure}
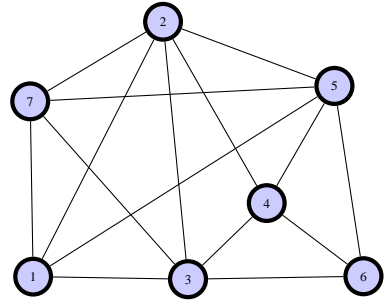

In this section, we evaluate the performance of the proposed diffusion algorithm through computer
simulations.
In all simulations, we consider the connected network of 7 agents shown in
Fig.\,\ref{fig:network} and employ the uniform combination rule $a_{\ell,k}=1/|\mathcal{N}_k|$~\cite{cattivelli2010diffusion}
and set the step-sizes across the agents to a uniform value, $\mu_k=\mu$.
To evaluate the MSD, we run $400$ experiments and compute the average MSD across these experiments for different approaches.

Some distributed techniques  that rely on EM techniques, of the form studied in \cite{nowak_EM_1212679,EM_dis_4558075,EM_clustering_5693298}, have been proposed for useful but different applications over sensor networks, such as clustering or density estimation, but not for missing data scenarios considered in this study.
We therefore compare our solution against  centralized processing techniques.
We also mentioned earlier that the missing data model considered in this paper is more general than the
models investigated in the missing data literature.
Therefore, to assess the performance of the proposed approach
against existing techniques, first we try to detect the missing position and then employ existing techniques that require knowledge of these positions.
Note that the suggested approach in this section may not be optimal. It may be possible to improve the performance of exiting approaches if we pursue \emph{joint} detection of the missing position and the estimation of the unknown parameter vector rather than follow the two-step procedure used in this section.
%
%

Suppose the $j-$th component of the regressor $\boldsymbol{u}_{k,i}$
at the $k-$th node at time $i$, denoted by $\boldsymbol{u}^j_{k,i}$,
is missing.
For such a scenario, $\boldsymbol{F}_{k,i}=\mathrm{diag}\{0,\ldots,\boldsymbol{f}^j_{k,i},\ldots,0\}$.
We define two hypotheses $\mathcal{H}_0$ and $\mathcal{H}_1$ as follows:
\begin{align}
\mathcal{H}_0: \bar{\boldsymbol{u}}^j_{k,i}=\boldsymbol{u}^j_{k,i}\\
\mathcal{H}_1: \bar{\boldsymbol{u}}^j_{k,i}=\boldsymbol{\xi}^j_{k,i}.
\end{align}
Since there are prior beliefs about $\mathcal{H}_0$ and $\mathcal{H}_1$, we consider the Bayesian hypothesis testing approach to detect the missing position~\cite{poor1994introduction}.
That is, we should decide $\mathcal{H}_1$ ($\boldsymbol{f}^j_{k,i}=1$ ) if the likelihood ratio (LLR) test is
larger than a threshold $\gamma$ as~\cite[Ch. 3]{kay1998fundamentals}:
\begin{align}
\label{eq:decesion_rule}
\text{LLR}\triangleq\frac{p_{\mathcal{H}_1}(\bar{\boldsymbol{u}}^j_{k,i})}{p_{\mathcal{H}_0}(\bar{\boldsymbol{u}}^j_{k,i})}\mathop{>}\limits^{\mathcal{H}_1} \gamma\triangleq\frac{(1-\widehat{p})}{\widehat{p}}
\end{align}
where $p_{\mathcal{H}_0}(\bar{\boldsymbol{u}}^j_{k,i})$ and $p_{\mathcal{H}_1}(\bar{\boldsymbol{u}}^j_{k,i})$ are the probability density
functions under hypotheses $\mathcal{H}_0$ and $\mathcal{H}_1$, respectively.

In general, evaluating the decision rule in \eqref{eq:decesion_rule} might be difficult.
In this example, we assume that $\boldsymbol{u}^j_{k,i}$ and $\boldsymbol{\xi}^j_{k,i}$ are normally distributed. Then,
\begin{align}
\label{eq:ga_dis}
p_{\mathcal{H}_0}(\bar{\boldsymbol{u}}^j_{k,i})&\,=\,\frac{1}{\sqrt{2\pi R_{u,k}(j,j)}}\exp\left[-\frac{|\bar{\boldsymbol{u}}^j_{k,i}|^2}{2R_{u,k}(j,j)}\right]\\
p_{\mathcal{H}_1}(\bar{\boldsymbol{u}}^j_{k,i})&\,=\,\frac{1}{\sqrt{2\pi \sigma^2_{\xi,k}}}\exp\left[-\frac{|\bar{\boldsymbol{u}}^j_{k,i}|^2}{2\sigma^2_{\xi,k}}\right].
\end{align}
The LLR is accordingly given by
\begin{align}
\text{LLR}=\sqrt{\frac{R_{u,k}(j,j)}{\sigma^2_{\xi,k}}}\exp{\left[\frac{|\bar{\boldsymbol{u}}^j_{k,i}|^2}{2R_{u,k}(j,j)}-\frac{|\bar{\boldsymbol{u}}^j_{k,i}|^2}{2\sigma^2_{\xi,k}}\right]}
\end{align}
and the decision rule can be expressed as follows:
\begin{align}
\label{eq:gau_rule}
|\bar{\boldsymbol{u}}^j_{k,i}|^2\left[\frac{1}{R_{u,k}(j,j)}-\frac{1}{\sigma^2_{\xi,k}}\right]\mathop{>}\limits^{\mathcal{H}_1}
 \log \left(\frac{(1-\widehat{p})^2\sigma^2_{\xi,k}}{\widehat{p}^2R_{u,k}(j,j)}\right).
\end{align}
It is still seen that for evaluating the decision rule \eqref{eq:gau_rule}, we need to know the variances $R_{u,k}(j,j)$ and $\sigma^2_{\xi,k}$.
To make the detection approach feasible, we assume that a good approximation for the ratio $R_{u,k}(j,j)/\sigma^2_{\xi,k}=r_k$ is available.
Note that for the proposed mATC algorithm, we do not need to know the ratio $r_k$.
From \eqref{eq:cov_rel}, we have
\begin{align}
\label{eq:var_est_cen_1}
(1-\widehat{p})R_{u,k}(j,j)+\widehat{p}\sigma^2_{\xi,k}=R_{\bar{u},k}(j,j)\approx \frac{1}{M_k}\sum_{i=1}^{M_k} \bar{\boldsymbol{u}}^j_{k,i}\bar{\boldsymbol{u}}^{*j}_{k,i}
\end{align}
where $M_k$ is the number of measurements collected by node $k$.
Substituting $R_{u,k}(j,j)$ by $r_k\sigma^2_{\xi,k}$ in \eqref{eq:var_est_cen_1}, we can
approximate the variance $\sigma^2_{\xi,k}$, denoted by $\check{\sigma}^2_{\xi,k}$, as
\begin{align}
\label{eq:var_estimate_cen_2}
\check{\sigma}^2_{\xi,k}=\frac{1}{M_k((1-\widehat{p})r_k+\widehat{p})}\sum_{i=1}^{M_k} \bar{\boldsymbol{u}}^j_{k,i}\bar{\boldsymbol{u}}^{*j}_{k,i}.
\end{align}
The variance $R_{u,k}(j,j)$ can be estimated as
\begin{align}
\label{eq:var_estimate_cen_3}
\widehat{R}_{u,k}(j,j)=r_k \check{\sigma}^2_{\xi,k}.
\end{align}
We note that the estimate of the variance $\sigma^2_{\xi,k}$ in $\eqref{eq:var_estimate_cen_2}$ needs the ratio $r_k$ to be known and also the estimate
is based on a batch processing.

Once the missing positions have been identified, we then apply two state-of-the-art centralized techniques to centrally estimate the unknown vector $w$.
In the first approach, the detected missing position is filled by the mean of the  data and then a least-squares construction is
applied (Imput-LS). For details on this approach, the reader may refer to \cite{little1992regression}.
In the second approach, we obtain the maximum likelihood estimator (MLE) assuming known distributions for
the regressor and perturbation $\boldsymbol{v}_{k}(i)$. The details of the MLE for missing data can be found in \cite{ibrahim2005missing,little1992regression,horton1999maximum,dempster1977maximum}.

\begin{remark}
The two-step approach above may not be an optimal implementation for the centralized solution. It may be possible to develop more efficient centralized algorithms based, for example, on the EM algorithm and mixture models \cite{bilmes1998gentle,Hastie_tibshirani,bishop:2006:PRML,EM_mixture,ghahramani1994supervised}. It is noted that if a mixture model for the missing data in \eqref{eq:relation_data} is considered, then a large number of components may be needed. 
\end{remark}

In the first simulation,
we assume a Gaussian distribution for the process noise, $\boldsymbol{v}_k(i)\sim \mathcal{N}(0,\sigma^2_{v,k})$ with $\sigma^2_{v,k}=0.01$.
The regressor $\boldsymbol{u}_{k,i}$ has Gaussian distribution with diagonal covariance matrix, $R_{u,k}=\mathrm{diag}\{1,1.6,0.8,0.95,1.2\}$.
The probability of missing is set to $p=0.3$ or $p=0.4$, which is assumed to be known in advance.
We set the unknown vector $w^o=[1,~-0.5,~1.2,~0.4,~1.5]^T$.
We also assume that $\boldsymbol{\xi}_{k,i}\sim\mathcal{N}(0,\sigma^2_{\xi,k})$ with
variances $\sigma^2_{\xi,1}= 0.02,$ $\sigma^2_{\xi,2}= 0.44,$ $\sigma^2_{\xi,3}=0.04,$ $\sigma^2_{\xi,4}=0.09$, $\sigma^2_{\xi,5}=0.15$,
$\sigma^2_{\xi,6}=0.26$, and $\sigma^2_{\xi,7}=0.13$. The step size $\mu$ is set to $0.04$ for every node.
In the simulation, we assume that the first component of the regressor is missing.
In this simulation, the recursions for estimating the variance $\sigma^2_{\xi,k}$, i.e., Eqs. \eqref{eq:var_es_1}--\eqref{eq:var_es_4}, start after $50$ iterations ($i=50$)
of Eqs. \eqref{eq:eq1}--\eqref{eq:eq2}.
We set $\alpha_1=\alpha_2=\alpha_3=0.01$.

\begin{figure}
\center
\subfigure[]{
\label{eq:MSD_1_p_3}
\psfrag{xlabel}[cc][][1]{Iteration}
 \psfrag{ylabel}[cc][][1]{MSD in dB}
 \psfrag{MLE}[cc][][.6]{MLE}
 \psfrag{mATC (Simulation)}[cc][][.6]{mATC (Simulation)}
 \psfrag{mATC (Theory)}[cc][][.6]{mATC (Theory)}
 \psfrag{LS-Imput}[cc][][.6]{Imput-LS}
  \includegraphics[width=80mm]{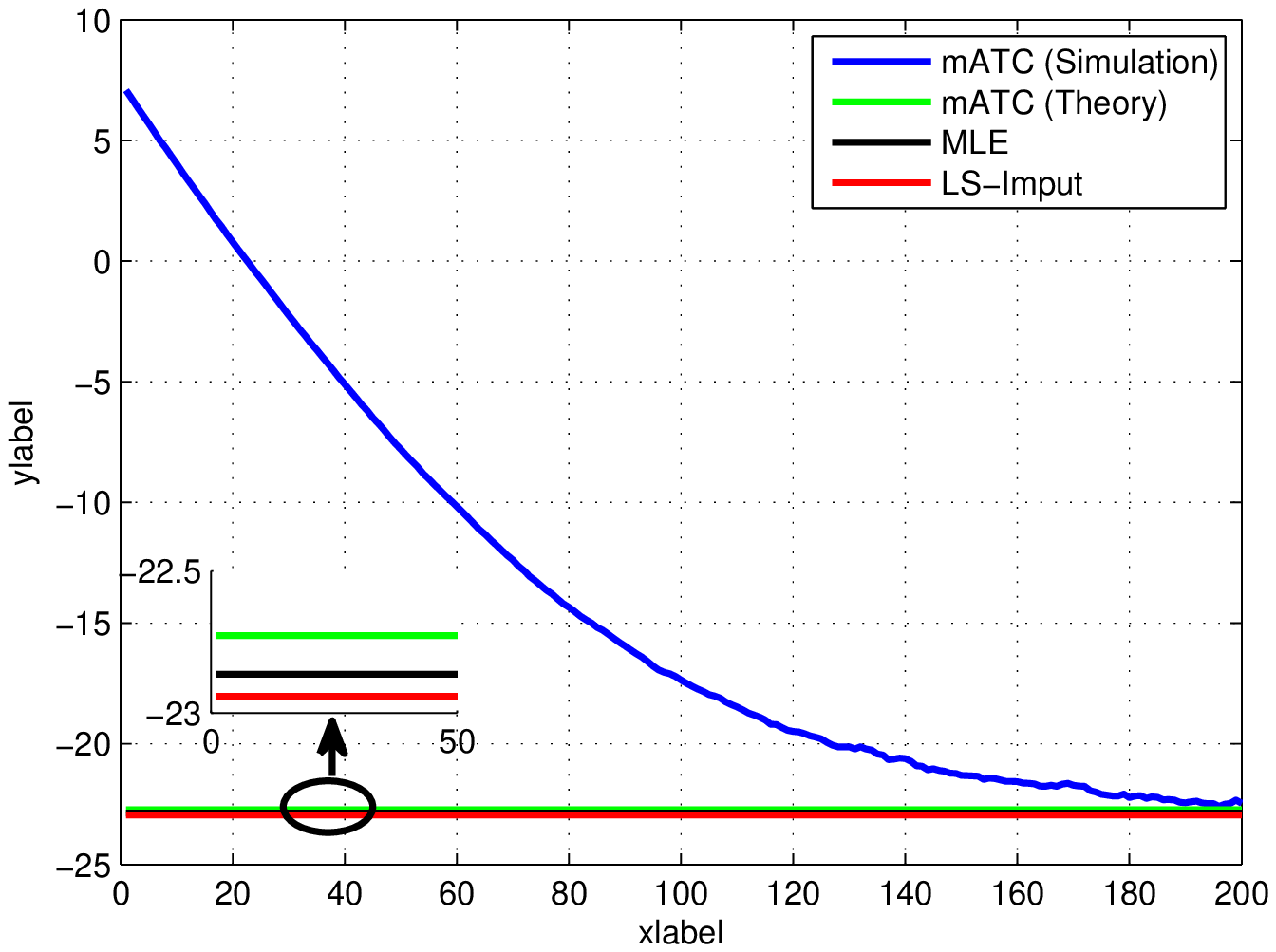}
  }
  \subfigure[]{
  \psfrag{xlabel}[cc][][1]{Iteration}
 \psfrag{ylabel}[cc][][1]{MSD in dB}
 \psfrag{MLE}[cc][][.6]{MLE}
 \psfrag{mATC (Simulation)}[cc][][.6]{mATC (Simulation)}
 \psfrag{mATC (Theory)}[cc][][.6]{mATC (Theory)}
 \psfrag{LS-Imput}[cc][][.6]{Imput-LS}
  \includegraphics[width=80mm]{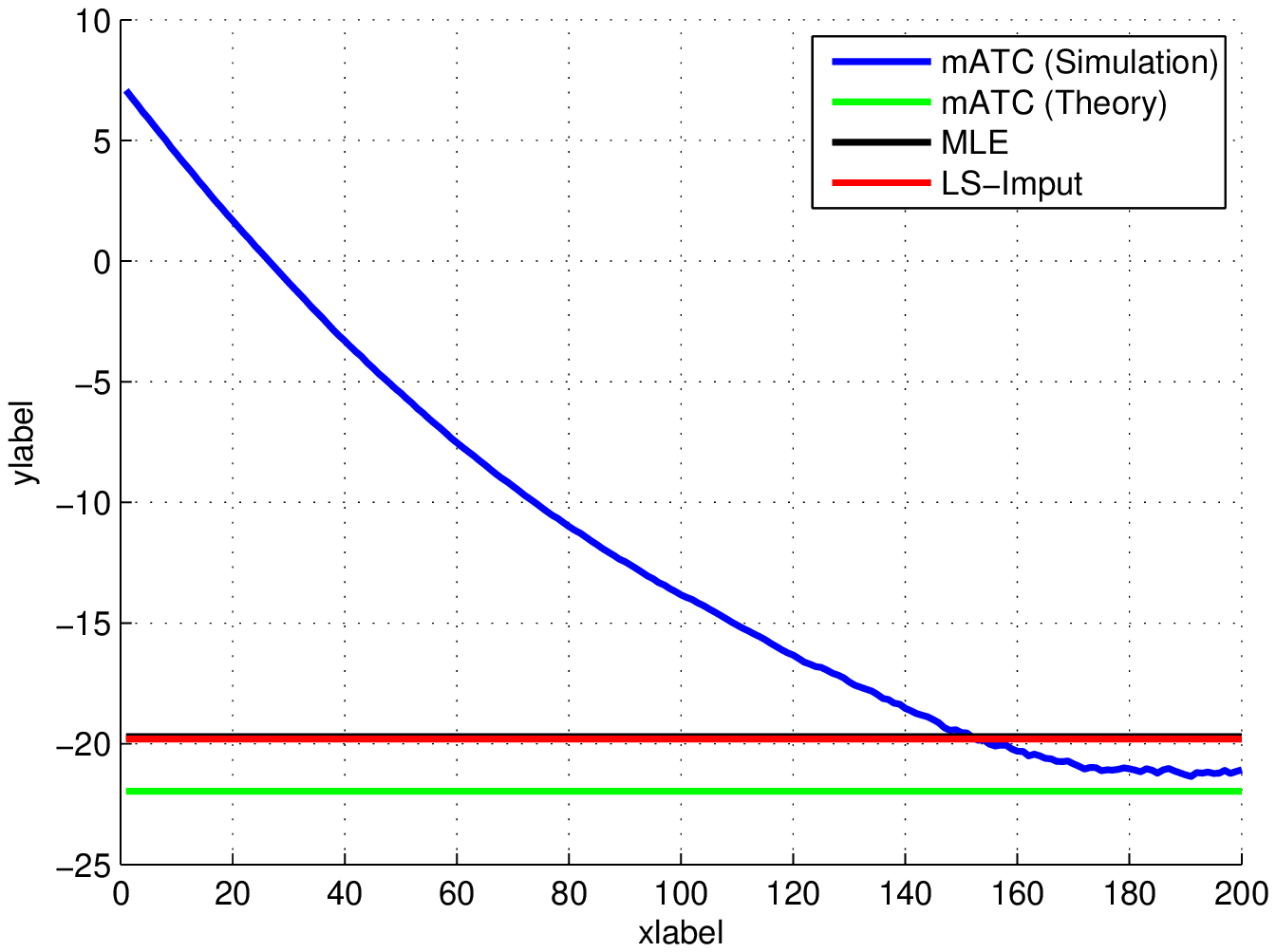}
  \label{eq:MSD_1_p_4}
  }
 \caption{The MSD learning curves for different diffusion algorithms at node $1$ for \subref{eq:MSD_1_p_3} $p=0.3$ and \subref{eq:MSD_1_p_4} $p=0.4$.}
 \label{fig:msd_different_scenarios}%
\end{figure}

Figures \ref{eq:MSD_1_p_3} and \ref{eq:MSD_1_p_4} show the MSD learning curves for different approaches for $p=0.3$ and $p=0.4$.
As it is observed, the mATC shows a promising performance compared to centralized approaches.
As the probability of missing increases, the proposed mATC considerably outperforms other approaches.
The reason is that the error of missing position detection will increase as the probability of missing increases.
Therefore, the performance of centralized approaches will be degraded more with increasing the probability of missing $p$.

In the next sections, we consider two specific applications.

\subsection{Household Consumption}
\label{eq:house_hould}
In this section, we evaluate the performance of the proposed algorithm for a household consumption application.
Household consumption depends on a number of parameters such as income, wealth, family size, and retirement status~\cite{rigobon2007estimation}.
It is assumed that the wealth variable is missed in the survey.
We consider the following log model for household consumption~\cite{rigobon2007estimation,parker2000spendthrift}:
\begin{align}
\label{eq:consump}
\ln \boldsymbol{c}_{k}(i)=&\alpha+(\ln \boldsymbol{l}_{k,i})\beta_1+(\ln \boldsymbol{m}^p_{k,i})\beta_2+(\ln \boldsymbol{m}^c_{k,i})\beta_3\nonumber\\
&+\boldsymbol{t}_{k,i}\beta_4+\boldsymbol{v}_{k}(i)
\end{align}
where $\alpha$ is $\ln \boldsymbol{c}_{k}(i)$ intercept, $\boldsymbol{c}_{k}(i)$ is the
consumption of household $k$ at time $i$, $\boldsymbol{l}_{k,i}$ is the total wealth, which is assumed to be censored, $\boldsymbol{m}^p_{k,i}$ is
the permanent part of the income
(a long-term measurement of average income that depends on a
number of parameters such as family income and education)
\cite{campbell1990permanent,lusardi1996permanent,parker2000spendthrift,hall1979stochastic,friedman1957},
$\boldsymbol{t}_{k,i}$ refers to the retirement status and family size. The modeling error $\boldsymbol{\epsilon}_{k}(i)$
is assumed to be zero-mean.
In a manner similar to \cite{rigobon2007estimation}, we only consider the first 4 components of the regressor, i.e., we set $\beta_4=0$.
As suggested by the earlier Remark \ref{re:rem1}, we subtract the mean of the measurements, which is given by
\begin{align*}
\mathbb{E}\,(\ln \boldsymbol{c}_{k}(i))=&\alpha+\mathbb{E}(\ln \boldsymbol{l}_{k,i})\beta_1+\mathbb{E}(\ln \boldsymbol{m}^p_{k,i})\beta_2\nonumber\\
&+\mathbb{E}(\ln \boldsymbol{m}^c_{k,i})\beta_3
+\mathbb{E}\boldsymbol{t}_{k,i}\beta_4
\end{align*}
from both sides of \eqref{eq:consump} and arrive at the model
\begin{align}
\boldsymbol{d}^c_{k}(i)
=\boldsymbol{u}^c_{k,i}w_c+\boldsymbol{v}_{k}(i)
\end{align}
where
\begin{align}
w_c&=[\beta_1~\beta_2~\beta_3]^T\nonumber\\
\boldsymbol{d}^c_{k}(i)&\triangleq \ln \boldsymbol{c}_{k}(i)-\mathbb{E}\,(\ln \boldsymbol{c}_{k}(i))\nonumber\\
\boldsymbol{u}^c_{k,i}&\triangleq [\ln \boldsymbol{l}_{k,i}~\ln \boldsymbol{m}^p_{k,i}~\ln \boldsymbol{m}^c_{k,i}]-\mathbb{E}\, [\ln \boldsymbol{l}_{k,i}~\ln \boldsymbol{m}^p_{k,i}~\ln \boldsymbol{m}^c_{k,i}].
\end{align}

Using a complete set of data, the authors in \cite{rigobon2007estimation} reported the
estimate ${\widehat{w}_c= [0.054,~0.182,~0.204]^T}$ for the unknown parameters.
We generate data according to $\widehat{w}_c$ and assume that the regressor $\boldsymbol{u}^c_{k,i}$ has Gaussian distribution. We model $\boldsymbol{v}_{k}(i)$ by a
zero-mean Gaussian random variable.
We further assume that the log of wealth is randomly missed and
we consider a uniform distribution over $[-\bar{q}, \bar{q}]$ for the missing variable, thus $\sigma_{\xi,k}^2=\bar{q}^2/3$.
In the simulation, we use $\bar{q}=0.5$.
In the
survey, it has been observed that approximately 30\%  of total wealth, including housing and stock market, is censored~\cite{rigobon2007estimation}, i.e., $p=0.3$.
In the simulation, we use $\mu=0.025$.
The updating step for estimating the variance $\sigma^2_{\xi,k}$ is executed from the beginning $i=1$.
We also set $\alpha_1=\alpha_2=0.001$, and $\alpha_3=0.0001$.

To derive a decision rule for this example, we first consider two distributions under hypotheses $\mathcal{H}_0$ and $\mathcal{H}_1$:
\begin{align}
\label{eq:ga_dis_uni}
p_{\mathcal{H}_0}(\bar{\boldsymbol{u}}^{cj}_{k,i})&\,=\,\frac{1}{\sqrt{2\pi R_{u^c,k}(j,j)}}\exp\left[-\frac{|\bar{\boldsymbol{u}}^{cj}_{k,i}|^2}{2R_{u^c,k}(j,j)}\right]\\
p_{\mathcal{H}_1}(\bar{\boldsymbol{u}}^{cj}_{k,i})&\,=
\left\{\begin{array}{ll}
\,\frac{1}{2{\bar{q}}},& \text{if}~ |\bar{\boldsymbol{u}}^{cj}_{k,i}|\leq {\bar{q}}\\
0,&\text{otherwise}
\end{array}
\right.
\end{align}

\noindent In the scenario for $|\bar{\boldsymbol{u}}^j_{k,i}|\leq {\bar{q}}$ the LLR is obtained as
\begin{align}
\label{eq:llr_2}
\text{LLR}=\frac{\sqrt{2\pi R_{u^c,k}(j,j)}}{2{\bar{q}}}\exp{\left[\frac{|\bar{\boldsymbol{u}}^{cj}_{k,i}|^2}{2R_{u^c,k}(j,j)}\right]}, \quad |\bar{\boldsymbol{u}}^j_{k,i}|\leq {\bar{q}}.
\end{align}

\noindent Since $\bar{q}$ and $R_{u^c,k}(j,j)$ are unknown \textit{a priori},
we need to approximate them from the data. Similar to the previous case,
we assume a good approximation of the ratio $r_k$ is available that helps us to obtain estimates of $\bar{q}$ and $R_{u^c,k}$ similar to \eqref{eq:var_estimate_cen_2}--\eqref{eq:var_estimate_cen_3}.
Then, we find an estimate for
$\widehat{\bar{q}}$ as
\begin{align}
\label{eq:uni_q_est}
\widehat{\bar{q}}&=\sqrt{3\check{\sigma}^2_{\xi,k}}\nonumber\\
&=\sqrt{\frac{3}{M_k\big((1-\widehat{p})r_k+\widehat{p}\big)}\sum_{i=1}^{M_k} \bar{\boldsymbol{u}}^{cj}_{k,i}\bar{\boldsymbol{u}}^{cj*}_{k,i}}.
\end{align}

\begin{figure}
\center
\psfrag{p1}[cc][][.7]{$(1-\widehat{p})\,p_{\mathcal{H}_0}(\bar{\boldsymbol{u}}^{cj}_{k,i})$}
 \psfrag{p2}[cc][][.7]{$\widehat{p}\,p_{\mathcal{H}_1}(\bar{\boldsymbol{u}}^{cj}_{k,i})$}
 \psfrag{h0}[cc][][.7]{$\mathcal{H}_0$}
 \psfrag{h1}[cc][][.7]{$\mathcal{H}_1$}
 \psfrag{q}[cc][][.7]{$\widehat{\bar{q}}$}
 \psfrag{-q}[cc][][.7]{$-\widehat{\bar{q}}$}
\includegraphics[width=85mm]{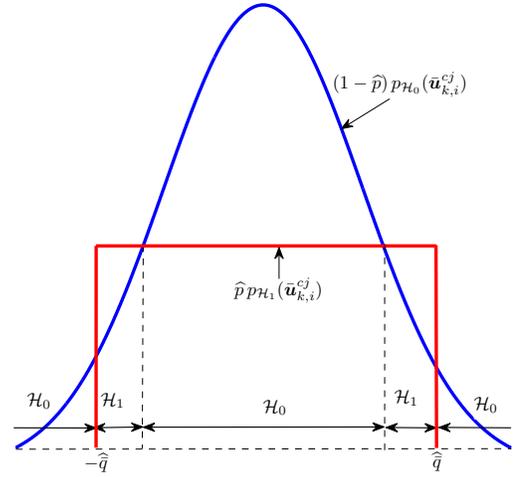}
 \caption{Decision regions and scaled probabilities with prior beliefs about $\mathcal{H}_0$ and $\mathcal{H}_1$.}
\label{fig:de_bound}%
\end{figure}

Using \eqref{eq:llr_2} and \eqref{eq:decesion_rule}, we obtain the following decision rule:
\begin{align}
\mathrm{if}~|\bar{\boldsymbol{u}}^j_{k,i}|\leq \widehat{\bar{q}}:
 \big|\bar{\boldsymbol{u}}^{cj}_{k,i}\big|^2\mathop{>
}\limits^{\mathcal{H}_1} 2 \widehat{R}_{u^c,k}(j,j)\log \left(\frac{2(1-\widehat{p})\widehat{\bar{q}}}{\widehat{p}\sqrt{2\pi \widehat{R}_{u^c,k}(j,j)}}\right)
\end{align}
where we replaced $R_{u^c,k}(j,j)$ and $\bar{q}$, respectively, by their estimates $\widehat{R}_{u^c,k}(j,j)$ and $\widehat{\bar{q}}$.
Figure \ref{fig:de_bound} shows an example of two distributions scaled by $\widehat{p}$ and $(1-\widehat{p})$ and the corresponding decision regions.

\begin{figure}
\center
\psfrag{xlabel}[cc][][1]{Iteration}
 \psfrag{ylabel}[cc][][1]{MSD in dB}
 \psfrag{MLE}[cc][][.6]{MLE}
 \psfrag{mATC (Simulation)}[cc][][.6]{mATC (Simulation)}
 \psfrag{mATC (Theory)}[cc][][.6]{mATC (Theory)}
 \psfrag{LS-Imput}[cc][][.6]{Imput-LS}
\includegraphics[width=80mm]{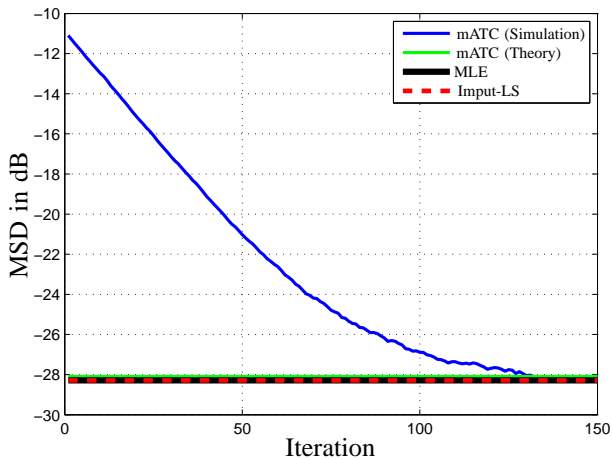}
 \caption{The MSD of different algorithms for the household data for node 1.}
\label{fig:msd_hosehold}%
\end{figure}

 Figure \ref{fig:msd_hosehold} shows the MSDs of the estimators for the household consumption data.
As it is observed, the proposed mATC shows comparable performance
with the centralized approach.

\subsection{Mental Health Survey}
\begin{table*}
\footnotesize
\centering
\caption{ Depression covariate for the $i$-th individual at the $k-$th company $\left(\widetilde{\boldsymbol{u}}_{k,i}=[\widetilde{\boldsymbol{u}}^1_{k,i},\widetilde{\boldsymbol{u}}^2_{k,i},\widetilde{\boldsymbol{u}}^3_{k,i},\widetilde{\boldsymbol{u}}_{k,i}^4,
\widetilde{\boldsymbol{u}}_{k,i}^5,\widetilde{\boldsymbol{u}}_{k,i}^6,\widetilde{\boldsymbol{u}}_{k,i}^7]\right)$ \cite{afifi2004computer}.}
 \begin{tabular}{| p{2.5cm}|p{15cm}|}
    \hline
    Covariate &   Range\\ \hline \hline
     Gender: $\widetilde{\boldsymbol{u}}^1_{k,i}$ &  =1 if the $i$-th responder is male;  =2 for the female\\
      Age in year: $\widetilde{\boldsymbol{u}}^2_{k,i}$& 19,20,\ldots,73\\
      Marital: $\widetilde{\boldsymbol{u}}^3_{k,i}$ & =1  never married; =2  married; =3  divorced;=4  separated; =5  widowed\\
     Education: $\widetilde{\boldsymbol{u}}^4_{k,i}$  & =1  less than high school; =2  some high
school; =3  finished high school; =4  some
college; =5  finished bachelor's degree;
=6  finished master's degree; =7  finished doctorate\\
       Log of Income: $\widetilde{\boldsymbol{u}}^5_{k,i}$ &   $\log 4000$ to $\log 55000$ \\
       Religion: $\widetilde{\boldsymbol{u}}^6_{k,i}$ &  =1 Protestant; =2  Catholic; =3  Jewish;
=4 none; =5  other\\
        Employment: $\widetilde{\boldsymbol{u}}^7_{k,i}$& =1  full time; =2  part time; =3  unemployed;
=4 retired; =5  houseperson; =6 in school;
=7 other\\
     \hline
    \end{tabular}
   \label{tab:regressor}
 \end{table*}
In the last simulation, we  consider the following model, motivated by a mental health survey study run by various companies\cite{afifi2004computer,little1992regression}:
\begin{align}
\label{eq:dep_model}
\bar{\boldsymbol{d}_k}(i)=\beta+\widetilde{\boldsymbol{u}}_{k,i}w+\boldsymbol{v}_k(i),\quad i=1,2,\ldots,N
\end{align}
where $\bar{\boldsymbol{d}_k}(i)$ is the square root of the total depression score
for every individual $i$, $\beta$ is the $\bar{\boldsymbol{d}_k}(i)$ intercept,
$\widetilde{\boldsymbol{u}}_{k,i}$ denotes the regressor (covariate) for every individual $i$, and $\boldsymbol{v}_k(i)$ is the modeling error.
Index $k$ refers to the company index and $i$ is used for participant's index.
The elements of $\widetilde{\boldsymbol{u}}_{k,i}$, defined in~\cite{afifi2004computer} are shown in Table \ref{tab:regressor} and they include variables such as income, age, and martial status.

The square root of the total depression score $\bar{\boldsymbol{d}_k}(i)$ is measured based on individual answers to 20
questions regarding feeling about depression~\cite{afifi2004computer}.
For example, the answer to the second question ``I felt depressed" can be a number between 0 and 3; 0 for less than 1 day feeling, 1 for one to two days, 2 for occasionally  or a moderate amount of time (3 to 4 days), and 3 for most of the time (5 to 7 days)~\cite{afifi2004computer}.

We apply the least squares technique to a subset of the data provided in \cite{afifi2004computer} to find an
estimate for $[\widehat{\alpha}~ \widehat{w}^T]$ in \eqref{eq:dep_model} as
$[\widehat{\alpha}~ \widehat{w}^T]=[ 0.1,0.27,-0.03,-0.06,0.13,0.73,-0.28,0.22]$.
We then use the estimate for $\widehat{w}$ to generate zero-mean square root total depression scores as follows.
Again, as indicated by the earlier Remark~1 and similar to the previous application, we modify the model of \eqref{eq:dep_model} as follows.
Consider the mean of both sides of model \eqref{eq:dep_model}, which is given by
\begin{align}
\label{eq:dep_model_2}
\mathbb{E} \bar{\boldsymbol{d}_k}(i)&=\beta_0+(\mathbb{E}\widetilde{\boldsymbol{u}}_{k,i}) w,
\end{align}
then, we subtract the above mean from both sides \eqref{eq:dep_model} to get
\begin{align}
\boldsymbol{d}_k(i)&\triangleq \bar{\boldsymbol{d}_k}(i)-\mathbb{E} \bar{\boldsymbol{d}_k}(i)\nonumber\\
&=\boldsymbol{u}_{k,i}w+\boldsymbol{v}_k(i),\quad i=1,2,\ldots,N,
\end{align}
where $\boldsymbol{u}_{k,i}=\widetilde{\boldsymbol{u}}_{k,i}-\mathbb{E} \widetilde{\boldsymbol{u}}_{k,i}$ is a zero-mean random vector.
To generate $\boldsymbol{u}_{k,i}$, we uniformly generate the regressor $\tilde{\boldsymbol{u}}_{k,i}$ according to
Table \ref{tab:regressor} and then subtract the mean.
We further assume that the income is missed with probability $0.3$ in the simulation study.
We consider a zero-mean Gaussian distribution with variance  $0.004$ for missing parts, i.e., $\sigma^2_{\xi,k}=0.004$.

The covariance matrix of the discrete regressor $\boldsymbol{u}_{k,i}$ is given by
\begin{align}
\label{eq:cov_mental}
R_{u,k}=\mathrm{diag}\{0.25,252,2,2.967,0.11,1.25,4\}.
\end{align}
The algorithm needs a smaller step size than the one in the previous simulation to converge due to
the largest eigenvalue of the covariance matrix in \eqref{eq:cov_mental}.
In the simulation, we set $\mu=0.0025$.
We assume a Gaussian distribution for measurement noise, $\boldsymbol{v}_k(i)\sim \mathcal{N}(0,\sigma^2_{v,k})$
with $\sigma^2_{v,k}=0.01$.
In this scenario, the updating for the estimate of
variance $\widehat{\boldsymbol{\sigma}}_{\xi,k}^2(i)$ starts from the beginning $i=1$.
Hence, the proposed algorithm is expected to have a slower convergence rate at smaller step-sizes. 
We set $\alpha_1=\alpha_2=\alpha_3=0.0001$.

\begin{figure}
\psfrag{xlabel}[cc][][1]{Iteration}
 \psfrag{ylabel}[cc][][1]{MSD in dB}
 \psfrag{MLE}[cc][][.6]{MLE}
 \psfrag{mATC (Simulation)}[cc][][.6]{mATC (Simulation)}
 \psfrag{mATC (theory)}[cc][][.6]{mATC (Theory)}
 \psfrag{LS-Imput}[cc][][.6]{Imput-LS}
 \includegraphics[width=80mm]{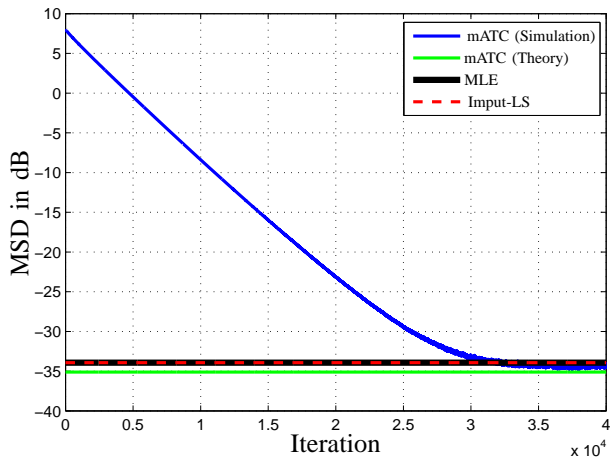}
\caption{The MSD of different approaches for the mental health survey data for node 1.}
 \label{fig:EMSE_ATC}
\end{figure}

Figure \ref{fig:EMSE_ATC} shows the MSDs of different approaches for mental health survey.
As it is observed the proposed mATC approach shows promising results compared to centralized techniques.
One way to improve the convergence rate is to estimate the curvature information of the cost function
and modify the update step \eqref{eq:eq1} based on the curvature information. 
Finally, Table \ref{tab:var_estimate} shows the estimate of the variance $\sigma^2_{\xi,k}$ at different nodes in steady state.
It is seen that every node can obtain a good estimate for $\sigma^2_{\xi,k}$ using the proposed mATC algorithm.

\begin{table}
\footnotesize
\centering
\caption{estimates of variance $\sigma^2_{\xi,k}$ in different nodes for the mental health survey data. The true value of the variance is $\sigma^2_{\xi,k}=0.004$. }
 \begin{tabular}{| l | l | l | l |l|l|l|}
   \hline
   Node 1& Node 2&Node 3&Node 4& Node 5&Node 6&Node 7\\ \hline\hline
   0.0047&  0.0046&    0.0048&    0.0049 &   0.0046&    0.0049&    0.0047\\\hline
 \end{tabular}
 \label{tab:var_estimate}
 \end{table}

\section{CONCLUSIONS}
\label{sec:conclude}
In this paper, we examined the estimation of an unknown vector
over a connected network of agents, with each agent subjected to a stream of data with incomplete regressors.
We have shown that the estimator in general is biased; hence, we have modified the cost function by  a (de)regularisation term
to mitigate the bias and obtained a distributed approach based on diffusion adaptation techniques.
We have also suggested a technique to estimate the (de)regularization term from the data.
We have studied the performance of the proposed algorithm under some simplifying assumptions and  considered two applications
in mental health and household consumption surveys.
Simulation results show a comparable performance
compared to existing centralized approaches based on imputation techniques.

\bibliographystyle{IEEEtran}
\bibliography{Ref}

\end{document}